\documentclass[twoside,english,11pt]{article}
\usepackage[T1]{fontenc}
\usepackage{color}
\usepackage{amsmath}
\usepackage{amssymb}
\usepackage{amsthm}
\usepackage{amsfonts}
\usepackage{graphicx}
\usepackage{grffile}
\usepackage{cancel}
\usepackage{verbatim}
\usepackage{multirow}
\usepackage{algorithm2e}[scright]
\usepackage[all]{xy}
\usepackage{hyperref}
\usepackage{multimedia}
\usepackage{multirow}
\usepackage{subfig}
\usepackage{placeins}
\usepackage{authblk}
\usepackage[T1]{fontenc}

\usepackage{calligra}

\usepackage{geometry}
\usepackage{comment}

\usepackage{pgfplots}
\usepackage{pgfplotstable}
\usepackage{tikz}
\usepgfplotslibrary{external}
\usepackage[colorinlistoftodos]{todonotes}
\makeatletter
\renewcommand{\todo}[2][]{\tikzexternaldisable\@todo[#1]{#2}\tikzexternalenable}
\makeatother

\usepackage{custom-commands-private}
\usepackage[squaren]{SIunits}

\newtheorem{problem}{Problem}

\SetKwComment{Comment}{$\triangleleft$ }{}

\makeatother

\title{Exploiting high-contrast Stokes preconditioners to efficiently
solve incompressible fluid--structure interaction problems}

\author[1]{Micha\l{} Wichrowski\footnote{Corresponding author: mt.wichrowsk@uw.edu.pl}}
\author[2]{Piotr Krzyżanowski}
\author[3]{Luca Heltai}
\author[4]{Stanis\l{}aw Stupkiewicz}

\affil[1]{Interdisciplinary Center for Scientifc Computing, Heidelberg University, Im Neuenheimer Feld 205, 69120 Heidelberg, Germany}
\affil[2]{Faculty of Mathematics, Informatics and Mechanics, University of Warsaw, Banacha 2, 02-097 Warsaw, Poland}
\affil[3] {Mathematics Area, SISSA – International School for Advanced Studies, via Bonomea 265, 34136, Trieste, Italy}
\affil[4]{Institute of Fundamental Technological Research, Polish Academy of Sciences, Pawińskiego 5B, 02-106 Warsaw, Poland}
\date{\vspace{-5ex}}

\providecommand{\keywords}[1]
{
  \small	
  \textbf{\textit{Keywords:}} #1
}

\usepackage{babel}
\begin{document}

\maketitle 

\begin{abstract}
In this work, we develop a new algorithm to solve large-scale incompressible
time-dependent fluid--structure interaction (FSI) problems using a matrix-free finite element method in
arbitrary Lagrangian--Eulerian (ALE) frame of reference. We derive a 
semi-implicit time integration scheme which improves the geometry-convective
explicit (GCE) scheme for problems involving the interaction between incompressible
hyperelastic solids and incompressible fluids. The proposed algorithm relies on
the reformulation of the time-discrete problem as a 
generalized Stokes problem with strongly variable coefficients, for which optimal
preconditioners have recently been developed. The resulting algorithm is
scalable, optimal, and robust: we test our implementation on model problems
that mimic classical Turek benchmarks in two and three dimensions, and
investigate timing and scalability results.
\end{abstract}
\keywords{
 fluid-structure interaction,
 finite element method,
 arbitrary Lagrangian-Eulerian,
 monolithic scheme,
 matrix-free method,
 geometric multigrid}

\section{Introduction}

The study of fluid-structure interaction (FSI) is crucial in numerous fields of
science and engineering. One of the earliest references to FSI is attributed
to Taylor~\cite{taylor1952fsi}, who described the interaction between a flexible
swimming animal and the surrounding fluid. Since then, FSI has remained an
active area of research, with significant advancements in computational methods \cite{fernandez2005newton,badia2008splitting,degroote2009performance,crosetto2011parallel},
time stepping techniques~\cite{crosetto2011fluid}, coupling strategies~\cite{benra2011comparison,bucelli2022fsi}, and preconditioning techniques~\cite{clevenger2019comparison,failer2020parallel}
(see, for example, the monographs by Bazilevs et al.~\cite{bazilevs2013computational} 
or Richter~\cite{richter2017fluid} for an overview on the topic).

The governing equations for FSI problems can be formulated in different
reference frames, such as the Lagrangian (classical for the solid), Eulerian
(classical for the fluid), or the arbitrary Lagrangian-Eulerian
(ALE)~\cite{donea2017arbitrary}. The Eulerian formulation is based on a fixed
coordinate system, where fluid and solid properties are computed at fixed
spatial positions. In the Lagrangian formulation, the equations of motion are
defined in a reference domain, and the domain deformation follows the motion of
the particles, while in the ALE formulation, the equations of motions are still
defined in a reference domain, but the deformation of the domain is not
necessarily attached to the motion of particles. All three approaches are
commonly used in solving FSI problems~\cite{richter2010finite,
richter2017fluid, ryzhakov2010monolithic, dunne2006eulerian}. While the ALE and
the Lagrangian formulations are \emph{body-fitted} methods, in the fully
Eulerian formulation, and in those formulations in which the fluid and the solid
equations are solved on separate domains, one needs a special treatment at the
interface, where the domains are usually \emph{non-matching}. In this case, one
needs to apply either fictitious domain methods~\cite{glowinski1994fictitious},
distributed Lagrange multiplier methods \cite{BoffiGastaldiHeltai-2018-a}, or
immersed methods~\cite{Peskin2002,BoffiGastaldiHeltai-2008-a}. In this paper, we
focus on the ALE formulation.

The combination of two (generally non-linear) continuum mechanics models
brings with it the inherent complexity of both of them, which is further
compounded by the challenges that arise due to their coupling. It is not
surprising that the solution of a fully-coupled monolithic FSI
problem~\cite{fernandez2005newton,hron2006monolithic,wick2017variational,richter2015monolithic}
is typically avoided by employing simplified coupling models, such as
partitioned loosely coupled schemes (explicit schemes)~\cite{burman2013explicit}
or partitioned fully coupled schemes (fixed-point or implicit
schemes)~\cite{causin2005added, kuttler2008fixed} (see, for example, references
\cite{benra2011comparison,landajuela2016,bucelli2022fsi} for a review of the
various methods).

In the case of FSI problems involving incompressible fluids, fully partitioned
schemes are often unstable due to problems with added-mass effect
\cite{causin2005added}. The natural remedy of leaning towards classical
fixed-point methods is often considered too expensive~\cite{gerbeau2005fluid},
however, some successful approaches have been developed~\cite{wall2008fluid}.

Implicit schemes lead to a system of non-linear equations to be solved at each
time step. The most complex computational scheme -- and most stable -- is the
fully implicit scheme~\cite{fernandez2005newton}, where the whole problem is
treated implicitly, and all Jacobians are computed exactly, resulting in a system
of equations that couples in a non-linear way the fluid motion, the deformation
of the solid, and the evolution of the domain geometry.
In~\cite{jodlbauer2019parallel, failer2020parallel} block preconditioners for
such a system were proposed.

Although this is indeed a very robust method, solving such a nonlinear problem
is not necessarily needed, since stable integration schemes can be obtained by
first using explicit methods to predict the deformed domain, and then use an
implicit scheme to compute the flow and deformation. This is the idea behind the
geometry-convective explicit (GCE) scheme
\cite{crosetto2011fluid,murea2017updated, badia2008modular, yang2015modeling}.
Well-posedness of the linear system arising from that discretization has been
proven in \cite{xu2015well}, and some block preconditioners have been developed. 

In this paper, we design and evaluate an efficient and scalable, fully-coupled, semi-implicit FSI solver
based on a stabilized finite element method, which uses parallel matrix-free
computing to deal with large problem sizes and exploits %
preconditioners
for high-contrast Stokes problems. %
The solver has much better stability properties than the GCE, while retaining a
relatively low cost as compared to fully implicit methods. It uses a
velocity-based formulation that allows for a natural coupling between the fluid
and solid equations. We model the fluid as an incompressible Newtonian fluid and
the solid as an incompressible Mooney--Rivlin hyperelastic material (see,
e.g.,~\cite{holzapfel2001biomechanics}) and complement the solid model with a
volume stabilization technique that improves incompressibility of the solid. 

One of the key ingredients of the solver is a semi-implicit scheme that
improves on the GCE scheme by applying similar ideas to the non-linear terms in
the fluid and solid equations, treating them semi-implicitly with a
predictor--corrector algorithm, crafted via a careful modification of backward
differentiation formulas (BDF). 
Such rewriting allows us to reinterpret the major step of the resulting improved GCE
scheme as a high-contrast Stokes problem, for which efficient %
preconditioners have recently been developed, e.g.~\cite{OlshanskiiPetersReusken2006,wichrowski2022,jodlbauer2022matrixfree}. %

Another advantage of the proposed FSI solver is its ability to make a consistent
use of matrix-free computing, a technique where matrix operations are performed
directly on the data, without the need to store the matrix explicitly,
thus increasing CPU cache efficiency and lowering the overall memory
footprint~\cite{Kronbichler2012}.
This is only possible because --- thanks to the connection with the generalized Stokes problem provided by our time-stepping method --- we are able to adapt to our problem a specialized preconditioner~\cite{wichrowski2022}, matrix-free by design, and suitable for high-contrast problems. %

The paper is organized as follows: in Section~\ref{sec:model}, we describe the
mathematical formulation of the FSI problem and in
Section~\ref{sec:discretization}  we introduce a semi-implicit
time integration scheme and the spatial discretization using the finite-element
method. Section~\ref{sec:preconditioner} discusses the preconditioned iterative
solvers employed to efficiently solve the resulting linear system of equations.
In Section~\ref{sec:numerical-results}, we present the numerical experiments
carried out to validate our method, including a test resembling Turek--Hron
\cite{turek2006proposal} benchmarks. We investigate the stability and efficiency of
our solver for a range of parameters, including mesh size and time step.
Finally, in Section~\ref{sec:conclusions}, we summarize the main findings of our
work and discuss some perspectives for future research.

\vfill

\section{Fluid--structure interaction model}
\label{sec:model}
\subsection{Weak form of the conservation equations}

Let us consider an initial (reference) configuration given by the domain $\Or\subset \Re^{d}$, consisting of two non-overlapping sub-domains: a fluid domain $\Orf\subset\Or$, and a solid domain
$\Ors\subset\Or$ ($\Or={\Ors\cup\Orf}$).
In FSI problems, one expects changes in both the solid configuration and the fluid domain in time. We denote
the actual domain at time $t\in[0,T]$ by $\Omega(t)$, with the convention that $\Omega(0)=\Or$. The deformed solid then occupies the domain $\Os(t)$, while the fluid domain $\Of(t)$ occupies the region $\Omega(t)\setminus\Os(t)$. When it will be clear from the context which instant $t$ is referred to, we will usually drop the time argument of the domains.

In a fixed frame of reference, we use plain fonts to indicate scalar Eulerian variables (e.g., $\pressure$) and indicate vector and tensor Eulerian variables using boldface characters (e.g., $\vc$ or $\cauchy$). 

We indicate global fields for both the solid and the fluid without subscripts, and indicate with subscript $\solid$ fields referring to the solid domain $\Os$, and with subscript $\fluid$ fields referring to the fluid domain $\Of$, i.e., 
\begin{equation}
\vc(t,\xc)=\begin{cases}
\vcs(t,\xc) & \xc\in\Os\\
\vcf(t,\xc) & \xc\in\Of.
\end{cases}
\end{equation}

The field $\vc$ represents the velocity of whatever type of particle happens to be at point $\xc$ at time $t$.
Similarly, we define a global Cauchy stress tensor $\cauchy$, and a  global density field for the coupled problem, given by
\begin{equation}
    \cauchy(t,\xc)=\begin{cases}
    \cauchys(t,\xc) & \xc\in\Os\\
    \cauchyf(t,\xc) & \xc\in\Of,
    \end{cases}
    \qquad
    \rho(t,\xc)=\begin{cases}
    \rhos(t,\xc) & \xc\in\Os\\
    \rhof(t,\xc) & \xc\in\Of,
    \end{cases}
\end{equation}
which will depend on the constitutive properties of the fluid and of the solid. For brevity, we will skip the time and space dependence of the variables when confusion is not possible.

By imposing conservation of mass and of momentum (see, e.g.,~\cite{richter2017fluid}), one obtains the following formal system of equations:
\begin{equation}
\begin{cases}
\rho\frac{\Df \vc}{\Df t}-\div  \cauchy  & =\body,\\
\frac{\Df \rho}{\partial t}+ \rho\div \vc & =0
\end{cases}
\qquad \text{ in } \Omega(t)\setminus \Gamma_i(t),
\label{eq:Cauchy-momentum}
\end{equation}
where $\body(t)$ is a given external force field (per unit volume), $\cauchy$
is the Cauchy stress tensor defined above, and $\Gamma_i(t)$ is the fluid--solid interface $\Gamma_i(t)=\partial\Os(t)\cap\partial\Of(t)$.
Conservation of angular momentum is guaranteed if the Cauchy stress tensor is symmetric. 

We indicate with ``$\grad$'' and ``$\div$'' the spatial differential operators with respect to the coordinates $\xc$ in the deformed configuration, and indicate with $\matDt{\bkappa}$ the material derivative, i.e., for a vector field $\bkappa$ and a scalar field $\alpha$, these are defined as time derivatives along the flow:
\begin{equation}
\matDt{\bkappa} =\Dt{}\bkappa+(\grad\bkappa)\vc, \qquad \matDt{\alpha} =\Dt{}\alpha+\vc\cdot\grad\alpha. 
\end{equation}

We partition the boundary $\Gamma(t)=\partial\Omega(t)$ into Neumann
$\Gamma_{N}(t)$ and Dirichlet $\Gamma_{D}(t)$ parts, and we set the following transmission and boundary conditions
\begin{equation}
\begin{cases}
    \cauchyf \ncf +\cauchys \ncs &= \mathbf{0} \text{ on }\Gamma_{i}\\
    \vcf &=\vcs \text{ on }\Gamma_{i}\\[.5cm]
    \vc & =\vc_{D}^{*}\mbox{ on }\Gamma_{D}\\
    \cauchy \nc  & = \traction^{*}\mbox{ on }\Gamma_{N}
\end{cases}
\label{eq:transmission-and-boundary-conditions}
\end{equation}
where $\vc_{D}^{*}(t)$ is a given prescribed velocity, $\traction^{*}(t)$ is a given prescribed traction, and $\nc$, $\ncs$, and $\ncf$ denote the unit outer normals to $\partial\Omega(t)$, $\partial\Os(t)$, and $\partial\Of(t)$, respectively. 

Using standard notations for Sobolev spaces, we indicate with $\spaceQ = L^2(\Omega)$,  $\spaceV_0 := \{\vc \in H^1(\Omega)^{d} \text{ s.t. } \vc = \mathbf{0} \text{ on } \Gamma_D\}$, and the affine space $\spaceVD := \{\vc \in H^1(\Omega)^{d} \text{ s.t. } \vc = \vc^*_D \text{ on } \Gamma_D\}$. With the conditions expressed in \eqref{eq:transmission-and-boundary-conditions}, for any time $t$, we can formally derive a global weak form of the conservation equations for both the solid and the fluid as

\begin{problem}[Weak form of the conservation equations]
    \label{pb:conservation-equations}
Given $\body(t) \in \spaceV_0'$ and $\traction^{*}(t) \in H^{-1/2}(\Gamma_N)^d$ for each time $t$ in $[0,T]$, find $(\vc, \rho) \in \spaceVD\times\spaceQ$ such that
\begin{equation}
\begin{aligned}
    &\int\limits _{\Omega}\rho\matDt{\vc}\cdot\testv\;\df \xc+\int\limits _{\Omega}\cauchy : \epsilon(\testv)\;\df \xc =\int\limits _{\Omega}\body\cdot\testv\;\df \xc+\int\limits _{\Gamma_{N}}\traction^{*}\cdot\testv\;\df s && \forall\testv\in \spaceV_0\\
    &\int\limits _{\Omega}\matDt{\rho} q\;\df \xc + \int\limits _{\Omega}\rho \div(\vc) q\;\df \xc = 0 && \forall q\in \spaceQ
\end{aligned}
\label{eq:weak-form-no-ale}
\end{equation}
where $\epsilon(\bkappa):=\frac12(\grad\bkappa+(\grad\bkappa)^T)$ denotes the symmetric gradient of a vector field $\bkappa$.
\end{problem}

In order to close the system, we will need to provide constitutive equations, discussed in Section~\ref{sec:constitutive}, initial conditions, and a suitable
representation for the evolution of the domain $\Omega(t)$, discussed in the next section.

\subsection{Arbitrary Lagrangian--Eulerian formulation}
\label{sec:ALE}

\begin{figure}
    \begin{center}
    \includegraphics[scale=0.6]{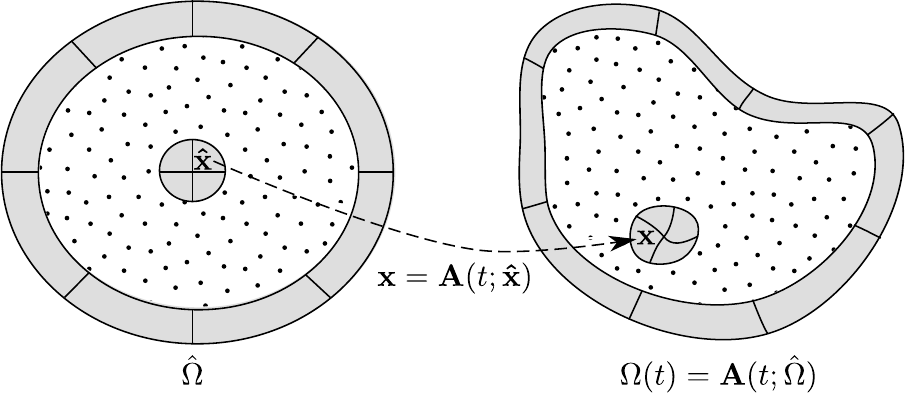} 
    \end{center}
    \caption{Initial domain $\Or$ and deformed domain $\Omega$(t). Solid
    marked with gray, fluid marked with dots, lines represent transformation
    of material points by $\mapas(t)$. \label{fig:phi-def}}
\end{figure}

We concretize the representation of the time-changing domain $\Omega(t)$, by introducing a diffeomorphism $\mapa(t):\Or\rightarrow\Re^{d}$ --- the arbitrary Lagrangian--Eulerian (ALE) map --- that maps points $\xr$ in the reference domain $\Or$ to points $\xc(t,\xr)=\mapa(t;\,\xr)$ on the deformed domain $\Omega(t)$, and such that $\Os(t) = \mapa(t;\Ors)$ and $\Of(t) = \mapa(t;\Orf)$. 

In general, we are free to choose the map $\mapa(t)$ arbitrarily, provided that the solid and fluid domains are mapped correctly at each time $t$. We choose an ALE map that coincides with the solid deformation map in the solid domain $\Ors$. On $\Orf$, we set $\mapa(t)$ as a pseudo-elastic extension in the fluid domain $\Orf$ (see, e.g.,~\cite{richter2017fluid}). Finally, we assume the ALE transformation equals the identity at time $t=0$, i.e., $\mapa(t=0; \xr) = \xr$ in $\Or\equiv\Omega(0)$.

We indicate with $\ur(t,\xr)=\mapa(t; \xr)-\xr$ the pseudo-displacement field that represents the domain deformation, and with $\wr$ its time derivative, representing the domain velocity and we adopt consistently the $\reference{\cdot}$ convention to indicate functions defined on $\Or$, or the pullback through $\mapa(t)$ of Eulerian fields defined on $\Omega(t)$ to $\Or$, i.e., $\reference{\bkappa} = \bkappa\circ\mapa$ and $\bkappa = \reference{\bkappa}\circ \mapa^{-1}$, and ``$\Grad$'' and ``$\Div$'' to indicate differential operators w.r.t.\ $\xr$. Subscripts $\solid,\fluid$ indicate the restriction of the field to the solid or fluid domains $\Ors,\Orf$ or $\Os(t), \Of(t)$.%

We use a Lagrangian setting in the solid domain, and let $\urs$ coincide with the displacement of solid particles in $\Ors$, so that $\partial_t\urs = \vrs$ in $\Ors$. In the domain $\Orf$, at every time $t$, we construct a pseudo-elastic extension $\urf$ of $\urs$ onto $\Orf$, namely 
\begin{equation}
\begin{split}\urf=\Ext(\urs)\end{split},\label{eq:arbitrary_Da}
\end{equation}
where $\Ext(\urs)$ is defined as the solution $\urf$ to
\begin{equation}
\begin{cases}
-\Div  (\mu_{A}\Symgrad(\urf)) & =\mathbf{0}\quad\mbox{ in }\Orf\\
\urf & = \urs\quad\mbox{on }\partial \Ors\\
\urf & = \mathbf{0} \quad\mbox{on }\partial \Orf\setminus\partial \Ors.
\end{cases}
\label{eq:extenstion}
\end{equation}
Here $\Symgrad(\urf) = \frac12(\Grad \urf + (\Grad \urf)^T)$ is the symmetric gradient of $\urf$ w.r.t.\ the coordinates $\xr$,  $\mu_{A}$ is a possibly variable coefficient, and for the sake of exposition we assume that the boundary of the fluid domain is not moving, even though arbitrary deformations could be applied to $\partial \Orf\setminus\partial \Ors$. Other (more computationally intensive) choices for the definition of $\Ext$ are possible, for example, using a biharmonic extension~\cite{helenbrook2003mesh,richter2017fluid}.

The equations of motion \eqref{eq:Cauchy-momentum} are complemented with an evolution equation for $\ur$ with zero initial conditions, defined through the domain velocity field $\wr$, i.e., 
\begin{equation}
    \label{eq:definition-domain-velocity}
    \wr := \begin{cases}
    \vr_s & \text{ in } \Ors\\    
    \Dt{\Ext(\urs)} & \text{ in } \Orf,    
    \end{cases}
    \qquad 
    \begin{cases}
    \Dt{\ur}(t,\xr) = \wr(t, \xr) & \text{ in } \Or\\
    \ur(t = 0,\xr) = \mathbf{0} & \xr\in\Or.
    \end{cases}
\end{equation}
Notice that, strictly speaking, only the equation for $\urs$ should be interpreted as an ODE with zero initial condition, while the equation for $\wrf$ is simply a time derivative evaluation.

We rewrite the weak formulation of the equations of motion in ALE form, by introducing the ALE time derivative in Eulerian coordinates, defined as the time derivative of a field at fixed $\xr$. For an Eulerian vector field $\bkappa$ or an Eulerian scalar field $\alpha$ we have
\begin{equation}
    \label{eq:ale-time-derivative}
    \begin{split}
    & \aleDt{\bkappa}(t, \mapa(t; \xr)) := \Dt{\reference{\bkappa}}(t, \xr) = \Dt{\bkappa}(t, \mapa(t; \xr)) + \left((\grad \bkappa)\wc\right)(t, \mapa(t; \xr)),\\
    & \aleDt{\alpha}(t, \mapa(t; \xr)) := \Dt{\reference{\alpha}}(t, \xr) = \Dt{\alpha}(t, \mapa(t; \xr)) + \left(\wc\cdot \grad \alpha\right)(t, \mapa(t; \xr)),
    \end{split}
\end{equation}
giving the following relation with the material time derivative:
\begin{equation}
    \label{eq:material-derivative-ale}
    \begin{split}
    & \matDt{\bkappa} = \aleDt{\bkappa} + (\grad\bkappa)(\vc-\wc), \\
    & \matDt{\alpha} = \aleDt{\alpha} + (\grad\alpha)\cdot(\vc-\wc),
    \end{split}
\end{equation}
which coincide with the partial time derivative in the Lagrangian case ($\vc=\wc$) and with the classical material time derivative in the purely Eulerian case ($\wc=0$).

We consider the mapped Sobolev spaces $\spaceQr := \spaceQ\circ\mapa$, $\spaceVr :=  \spaceV\circ\mapa$, and define, moreover, the space $\spaceVr_0 := H^1_0(\Or)^d$.

The weak form of the complete fluid-structure interaction problem in ALE form then reads
\begin{problem}[ALE weak form of conservation equations]
    \label{pb:ale-conservation-equations}
    Given $\body \in \spaceV_0'$ and $\traction^{*} \in H^{-1/2}(\Gamma_N)^d$, for each time $t$ in $[0,T]$, find $(\vc, \rho, \ur) \in \spaceVD\times\spaceQ\times\spaceVr$ such that
    \begin{equation}
        \begin{aligned}
            &\int\limits _{\Omega}\rho\left(\aleDt{\vc}+(\grad\vc)(\vc-\wc)\right)\cdot\testv\;\df \xc+\int\limits _{\Omega}\cauchy :\epsilon(\testv)\;\df \xc =\int\limits _{\Omega}\body\cdot\testv\;\df \xc+\int\limits _{\Gamma_{N}}\traction^{*}\cdot\testv\;\df s && \forall\testv\in \spaceV_0\\
            &\int\limits _{\Omega}\left(\aleDt{\rho}+(\vc-\wc)\cdot\grad \rho\right) q\;\df \xc + \int\limits _{\Omega}\rho \div(\vc) q\;\df \xc = 0 && \forall q\in \spaceQ\\
            &\int\limits _{\Orf} \mu_A \Symgrad (\ur ):\Symgrad (\testvr) +  \int\limits _{\Ors} \left(\Dt{\ur}-\wr\right) \cdot\testvr\;\df\xr = 0 && \forall \testvr\in \spaceVr,
        \end{aligned}
        \label{eq:weak-form-ale}
    \end{equation}
    with $\wc$ defined in \eqref{eq:definition-domain-velocity}.

    For simplicity of exposition, the first two equations in~\eqref{eq:weak-form-ale} are written in Eulerian form, but are solved on the reference domain $\Or$, using the spaces $\spaceVr$, and $\spaceQr$.
\end{problem}

From the numerical point of view, it may be convenient to express the dependency between $\ur$ and $\wr$ in a strong form, rather than using the weak form presented in~\eqref{eq:weak-form-ale}. Time discretization and how to write the time evolution of the grid is presented in Section~\ref{sec:Time-derivative-approximation}.

\subsection{Constitutive models}
\label{sec:constitutive}

\subsubsection{Incompressible Newtonian fluid}

We consider a classic incompressible Newtonian fluid with constant density $\rhof$ so that 
\begin{equation}
\cauchyf =2\eta_{\fluid}\epsilon(\vcf)-\pressuref \Id,\label{eq:const_fluid}
\end{equation}
where $\pressuref=-\frac13\trace\cauchyf$ is the fluid pressure, $\eta_{\fluid}$ is the viscosity and $\Id$ is the identity tensor. 
In this particular case, mass conservation given in equation~\eqref{eq:Cauchy-momentum} reduces to volume conservation, viz.
\begin{equation}
\div \vcf=0.\label{eq:fluid-divergence}
\end{equation}

We note that the pressure $\pressuref$ in Eq.~\eqref{eq:const_fluid} plays the role of a Lagrange multiplier enforcing the incompressibility constraint \eqref{eq:fluid-divergence}, and is solved for in place of the density $\rhof$ which is, in this constitutive model, a material constant of the fluid.

\subsubsection{Incompressible Mooney--Rivlin solid}
\label{sec:MR-solid}

We assume that the solid is an incompressible hyperelastic material governed by the Mooney--Rivlin model~\cite{ogden1972large,bigoni2012nonlinear}.
In the solid domain we define the deformation gradient 
\begin{equation}
\Fr :=\Grad\mapas=\Id+\Grad\urs,\label{eq:hat_f_def}
\end{equation}
and we have $\Fc = \Fr\circ\mapas^{-1}$.
The inverse deformation gradient is 
\begin{equation}
\Fc^{-1}=\grad\left(\mapas^{-1}\right)=\Id-\grad \ucs. \label{eq:f_inv}
\end{equation}

Incompressibility implies the following constraint on the determinant of $\Fr$, $\Jr = \det(\Fr)$ ($J=\Jr\circ\mapas^{-1}$),
\begin{equation}
\Jr = 1 ,
\end{equation}
which is equivalent to  $\div\vcs=0$.

In the case of an incompressible Mooney--Rivlin solid, the Cauchy stress can be expressed using the left Cauchy--Green deformation tensor $\Bc:=\Fc\Fc^{T}$ as follows~\cite{bigoni2012nonlinear}
\begin{equation}
\cauchys=\mu_{1}\Bc-\mu_{2}\Bc^{-1}-\pressures^{*}\Id,\label{eq:Cauch-mooney}
\end{equation}
where $\pressures^{*}$ is again a Lagrange multiplier that enforces incompressibility and $\mu_{1}\geq0$ and $\mu_{2}\geq0$ are material parameters such that $\mu_{\solid}=\mu_{1}+\mu_{2}>0$ is the shear modulus in the reference configuration. 
For $\mu_{2}=0$, the incompressible neo-Hookean model is obtained as a special case of the incompressible Mooney--Rivlin model.

The Cauchy stress~\eqref{eq:Cauch-mooney} can be equivalently expressed as
\begin{equation}
\cauchys=\mu_{1}(\Bc-\Id)+\mu_{2}(\Id-\Bc^{-1})-\pressures\Id, \qquad
\pressures=\pressures^{*}-\mu_{1}+\mu_{2},
\label{eq:Cauch-shifted}
\end{equation}
where $\pressures$ is a shifted Lagrange multiplier that vanishes whenever $\cauchys=\boldsymbol{0}$. 
Note, however, that $\pressures$ is not equal to the solid pressure, i.e., $\pressures\neq-\frac13\trace\cauchys$.

Substituting Eqs.~\eqref{eq:hat_f_def} and~\eqref{eq:f_inv} into Eq.~(\ref{eq:Cauch-shifted}) the Cauchy stress over the solid domain can now be expressed in the following form:
\begin{multline}
\cauchys = \mu_{1}\left(\left(2\Symgrad (\urs)+\Grad\urs(\Grad\urs)^T\right)\circ\mapas^{-1}\right)\\
+\mu_{2}\left(2\epsilon(\ucs)-(\grad \ucs)^{T}\grad \ucs\right) +
\pressures \Id .
\label{eq:solid-stress}
\end{multline}
We note that the second and the third term on the right-hand side are evaluated entirely in the current configuration. 
At the same time, the first term involves evaluation of the displacement gradient in the reference configuration, and the corresponding terms are then pushed forward through $\mapas^{-1}$ to the current configuration, where the integration is performed. 
In the present implementation, we restrict ourselves to the special case of $\mu_{1}=0$ and $\mu_{2}=\mu_{\solid}$, so that the first term in Eq.~\eqref{eq:solid-stress} vanishes. 
Note that in 2D plane-strain problems the incompressible Mooney--Rivlin model is equivalent to the neo-Hookean model with one shear modulus $\mu_{\solid}=\mu_{1}+\mu_{2}$ \cite{bigoni2012nonlinear}, hence the above assumption affects only 3D cases.

\subsection{Incompressible fluid-structure interaction problem in ALE form}

Let us summarize the results from this section by presenting the weak
form of the FSI problem in the ALE frame of reference. In the incompressible formulation we are using in this work, the conservation of mass transforms to conservation of volume. The (possibly different) densities $\rhos$ and $\rhof$ are constitutive constants of the fluid and of the solid, and the primal unknowns are the velocity field $\vc$, the pressure-like field $\pressure$ (acting as a Lagrange multiplier for the incompressibility constraints on both the solid and the fluid) and the pseudo-displacement $\ur$, corresponding to the solid displacement in $\Os$ and to the domain displacement in $\Of$.

\begin{problem}[ALE weak form of the incompressible fluid--structure interaction problem]
    \label{pb:ale-fsi-problem}
    Given $\body \in \spaceV_0'$ and $\traction^{*} \in H^{-1/2}(\Gamma_N)^d$, for each time $t$ in $[0,T]$, find $(\vc, \pressure, \ur) \in \spaceVD\times\spaceQ\times\spaceVr$ such that
    \begin{equation}
        \begin{aligned}
            &\int\limits _{\Omega}\rho\left(\aleDt{\vc}+(\grad\vc)(\vc-\wc)\right)\cdot\testv\;\df \xc+\int\limits _{\Omega}\cauchy :\epsilon(\testv)\;\df \xc  = && \\
            & \hspace{4cm}\int\limits _{\Omega}\body\cdot\testv\;\df \xc+\int\limits _{\Gamma_{N}}\traction^{*}\cdot\testv\;\df s && \forall\testv\in \spaceV_0\\
            &\int\limits _{\Omega}\div(\vc) q\;\df \xc = 0 && \forall q\in \spaceQ\\
            &\int\limits _{\Orf} \mu_A \Symgrad (\ur ):\Symgrad (\testvr) +  \int\limits _{\Ors} \left(\Dt{\ur}-\vr\right) \cdot\testvr\;\df\xr = 0 && \forall \testvr\in \spaceVr,
        \end{aligned}
        \label{eq:weak-form-fsi-ale}
    \end{equation}
    where $\wrf = \Dt{\Ext(\urs)}$, $\wcs = \vcs$, and the Cauchy stress $\cauchy$ takes the form:
    \begin{equation}
        \label{eq:stress-ale}
        \cauchy = \cauchy(\vc,\uc,\pressure) = \pressure \Id +  
        \begin{cases}
            2\eta_{\fluid} \symgrad(\vcf) & \text{ in } \Of\\
            2\mu_{\solid} \symgrad(\ucs) - \mu_{\solid} (\grad \ucs)^T\grad \ucs & \text{ in } \Os.
        \end{cases}
    \end{equation}
    The subdomains $\Os$ and $\Of$ appearing in integrals
    in Eq.~\eqref{eq:weak-form-fsi-ale} are defined as the images of the undeformed
    ones: $\Os=\mapa(\Ors)$ and $\Of=\mapa(\Orf)$, even though the computation is formally performed on the reference domain $\Or$.
    Additionally, the following initial conditions have to be met: 
    \begin{equation}
    \begin{cases}
    \vcs(t=0, \, \xr)  =\vcs^{*}&\text{ in }\Ors\\
    \urs(t=0, \, \xr) =\ucs^{*}&\text{ in }\Ors\\
    \vcf(t=0, \, \xr) =\vcf^{*}& \text{ in }\Orf.
    \end{cases}\label{eq:boundary_cond-1-1-1}
    \end{equation}
\end{problem}

This formulation is similar to the ones appearing in the literature
\cite{murea2017updated,bazilevs2013computational,xu2015well}, the major difference being the incompressibility assumption on the solid constitutive equations. At this stage the problem is still fully non-linear. The most difficult non-linearity comes from the fact that all domains are dependent on the (unknown) mapping $\mapa$, which in turn is hidden also in the definition of the ``$\grad$'' and ``$\div$'' differential operators, when they are expressed in the reference domain $\Or$.

\section{Discretization 
}
\label{sec:discretization} 

\subsection{Time integration scheme}
\label{sec:Time-derivative-approximation} 

We consider a uniform time discretization of the time
interval $[0,T]$ with a set of equidistant points $\{t^{0},...,t^{N}\}$
where $t^{n}=n\timestep$, with the time step size $\timestep=T/N$.
At the $n$-th time step $t^{n}$, we define $\vr^{n}$, $\urs^{n}$
and $\ur^{n}$ as an approximation of $\vr$, $\urs$
and $\ur$, respectively. The approximate  pseudo-displacement defines the approximation $\mapa^{n}$ of the mapping $\mapa(t^n; \cdot)$:
\begin{equation}
\label{eq:mapa-approx}
\mapa(t^{n};\xr)\approx \mapa^{n}(\xr)=\xr+\ur^{n}(\xr)\qquad\xr\in\Or,
\end{equation}
which in turn defines the approximate domains:
\begin{equation}
\label{eq:domain-approx}
\Omega(t^{n})\approx\Omega^{n}:=\mapa^{n}(\Or), \qquad
\Os(t^{n})\approx\Os^{n}:=\mapa^{n}(\Ors), \qquad 
\Of(t^{n}) := \Omega^{n}\setminus \Os^n.
\end{equation}
At the time step $n$, the relation between the point $\xc^{n}\in\Omega^{n}$
and the point $\xr\in\Or$ is defined by the mapping $\mapa^{n}$,
i.e: $\xc^{n}(\xr)=\mapa^{n}(\xr).$ That is, the spatial derivatives
\begin{equation}
\grad\kappa=\frac{\partial}{\partial x}\kappa,\qquad\epsilon(\kappa)=\frac{1}{2}\left(\grad\kappa+\left(\grad\kappa\right)^{T}\right)
\end{equation}
are associated with transformation $\mapa^{n}$. Since $\mapa^{n}$ is defined
by a pseudo-displacement, those derivatives are implicitly defined by $\ur^{n}$. 

We approximate the time derivative of an arbitrary field $\hat{\kappa}$ defined on $\Or\times[0,T]$ by the $k$-step backward differentiation formula (BDF) \cite{cash1980integration},
\begin{equation}
\Dt{\reference{\kappa}(t=t^n,\xr)}\approx\Dtau{\reference{\kappa}^{n}}(\xr) := \frac{1}{\gamma\timestep}\sum\limits _{i=0}^{k}\alpha_{i}\reference{\kappa}^{n-i}(\xr),
\label{eq:stencil}
\end{equation}
with the coefficients normalized so that $\alpha_{0}=1$. Here, we restrict ourselves to $k\leq 2$. The case $k=1$ corresponds to the classical implicit Euler scheme ($\gamma=1$, $\alpha_{1}=-1$), while for $k=2$ there holds $\gamma=\frac{2}{3}$, $\alpha_{1}=-\frac{4}{3}$, $\alpha_{2}=\frac{1}{3}$.

Therefore, at $t=t^n$, we will have
$$
\Dt{\ur}\approx\Dtau{\ur^{n}}, \qquad  \wrf = \Dt{\Ext(\urs)}\approx\wrf^n := \Dtau{\Ext(\urs^{n})}.
$$

For the ALE time derivative of the velocity we follow the identity $\aleDt{\vc}(t,\xc) = \Dt{\vr}(t,\xr)$ and approximate, cf. e.g. \cite{crosetto2011parallel,xu2015well,richter2017fluid}, at $\xc^n = \mapa^n(\xr)$, by moving back to the reference configuration:
\begin{equation}
\aleDt{\vc}(t=t^n,\xc^n) \approx\aleDtau{\vc^{n}}(\xc^n) := \Dtau{\vr^{n}}(\xr) = \frac{1}{\gamma\timestep}\sum_{i=0}^{k}\alpha_{i}\vc^{n,i}(\xc^n),
\label{eq:def-delta-v}
\end{equation}
where we set $\vc^{n,i}:=\vc^{n-i}\circ \mapa^{n-i}\circ(\mapa^{n})^{-1}$. %

After replacing the time derivatives %
in \eqref{eq:weak-form-fsi-ale} by their approximates,
we obtain an implicit time integration scheme:

\begin{problem}[Fully implicit BDF time discretization of ALE weak form of the incompressible fluid--structure interaction problem]
\label{pb:fully-implicit-ale-fsi-problem}
Given $\body \in \spaceV_0'$ and $\traction^{*} \in H^{-1/2}(\Gamma_N)^d$, for each time $t^n$, $n=1,\ldots,N$, find $(\vc^n, \pressure^n, \ur^n) \in \spaceVD^n\times\spaceQ^n\times\spaceVr$ such that
\begin{equation}
\begin{aligned}
&\int\limits _{\Omega^n}\rho\left(\aleDtau{\vc^n}+(\grad\vc^n)(\vc^n-\wc^n)\right)\cdot\testv\;\df \xc+\int\limits _{\Omega^n}
\cauchy^n :\epsilon(\testv)\;\df \xc = && \\
& \hspace{4cm}\int\limits _{\Omega^n}\body\cdot\testv\;\df \xc+\int\limits _{\Gamma_{N}}\traction^{*}\cdot\testv\;\df s && \forall\testv\in \spaceV_0^n\\
&\int\limits _{\Omega^n}\div(\vc^n) q\;\df \xc = 0 && \forall q\in \spaceQ^n\\
&\int\limits _{\Orf} \mu_A \Symgrad (\ur^n ):\Symgrad (\testvr) +  \int\limits _{\Ors} \left(\Dtau{\ur^n}-\vr^n\right) \cdot\testvr\;\df\xr = 0 && \forall \testvr\in \spaceVr,
\end{aligned}
\label{eq:fully-implicit-weak-form-fsi-ale}
\end{equation}
with $\wrf^n = \Dtau{\Ext(\urs^{n})}$ and $\wcs^n = \vcs^n$; the stress $\cauchy^n=\cauchy(\vc^n,\uc^n,\pressure^n) $ specified in \eqref{eq:stress-ale} and  the initial data prescribed by \eqref{eq:boundary_cond-1-1-1}. Here $\spaceVD^n$ and $\spaceQ^n$ denote the corresponding spaces over domain $\Omega^n$, which is implicitly defined by \eqref{eq:mapa-approx}--\eqref{eq:domain-approx}.
\end{problem}
Note that the last equation of \eqref{eq:fully-implicit-weak-form-fsi-ale} effectively splits into two,
$$
\urf^{n} =\Ext(\urs^{n}), \qquad \Dtau{\urs^{n}} =\vrs^{n}.
$$
From the latter and the linearity of $\Ext$ it follows that \begin{equation}
\wrf^n = \Dtau{\Ext(\urs^{n})} =  \Ext(\Dtau{\urs^{n}}) = \Ext(\vrs^n).
\end{equation}

\subsubsection{Semi-implicit scheme}
Problem~\ref{pb:fully-implicit-ale-fsi-problem}, which needs be solved on every time step, is nonlinear and usually is solved by the Newton's method (\cite{richter2015monolithic,heil2008solvers,langer2018numerical}) or a fixed-point method (without significant impact on stability, as demonstrated in \cite{lozovskiy2015unconditionally}).

Here, we take a different approach inspired by \cite{crosetto2011parallel,xu2015well,badia2008modular}, and simplify the nonlinear problem by splitting the solid displacement into an explictly predicted displacement and an implicit velocity dependent part. We exploit the structure of 
Problem~\ref{pb:fully-implicit-ale-fsi-problem} 
and in particular of the Cauchy stress definition~\eqref{eq:stress-ale} to expose a solution strategy based on fixed point iterations defined through a semi-implicit splitting of the BDF scheme defined in \eqref{eq:stencil} for a simpler problem.
The current pseudo-displacement and the current ALE velocity are affinely equivalent, and we expose this dependency together with the dependency on previously computed solutions ($\ure, \wre$) by expressing the currently unknown velocity and pseudo displacement $(\uri, \wri)$ as
\begin{equation}
  \label{eq:velocity-splitting}
  \uri := \gamma\timestep \wri - \sum\limits _{i=1}^{k}\alpha_{i}\ur^{n-i} =  \gamma\timestep \wri + (\ure - \gamma\timestep \wre)
\end{equation}
where $\wc^{\#}$ is the most recently computed velocity (i.e., the velocity from the previous time step or from a previous fixed point iteration) and 
 \begin{equation}
 \ure :=  \gamma\timestep \wre - \sum\limits _{i=1}^{k}\alpha_{i}\ur^{n-i}
 \end{equation}
 represents an explicit approximation of the pseudo-displacement.

This splitting is based on an explicit prediction $\uce$ of the next displacement, to be used in the computation of a temporary $\mapae$, and thus $\Oe$, and in an explicit computation of non-linear correction terms of the solid Cauchy stress. The remaining velocity dependent part is solved for implictly, resulting in the following splitting the global Cauchy stress
\begin{equation}
  \label{eq:stress-ale-semi-implicit}
  \cauchy = \pressure \Id +  2\eta \symgrad(\vci) + 
  \begin{cases}
      0 & \text{ in } {\explicit{\Of}} \\
      \mu_{\solid} \left[ 2 \symgrad\left(\sum\limits _{i=1}^{k}\alpha_{i}\uc^{n,i}\right) -(\grad \ucse)^T\grad \ucse\right] & \text{ in } {\explicit{\Os}},
  \end{cases}
\end{equation}
where 
\begin{equation}
  \label{eq:viscosity-splitting}
  \eta = 
  \begin{cases}
      \eta_f & \text{ in } {\explicit{\Of}}\\
      \gamma \timestep \mu_{\solid} & \text{ in } {\explicit{\Os}},
  \end{cases}
\end{equation} 
making the problem equivalent to a Stokes-like system on the domain $\explicit{\Omega}$, with jumps in the viscosities across the interface $\Gamma_i$. By explicitly writing the system in this way, we can exploit the robust preconditioner for Stokes-like systems with discontinuous viscosities developed in~\cite{wichrowski2022} also for the solution of FSI problems; cf.\ Section~\ref{sec:preconditioner}.

The final scheme, \GCSI{k}{M_1,\ldots,M_\nofmodes} (\textit{Geometry-Convective semi-Implicit}), which computes the solution on the next time step, is defined in Algorithm~\ref{alg:1}; its concept \emph{loosely} resembles an $\nofmodes$--stage predictor--corrector scheme.  On the $s$-th stage, $s=1,\ldots,\nofmodes$, the algorithm proceeds according to prescribed operational mode $\mode_s$. We consider two modes: if $\mode_s = \predictor$, the scheme computes the geometry and convective terms in an explicit way, while for  $\mode_s = \corrector$, it treats the convective term more implicitly (the details are provided later in this Section). The simplest and cheapest scheme of this kind, \GCSI{k}{\predictor}, coincides with the $k$-th order GCE scheme of \cite{crosetto2011parallel}, so our scheme can be considered a generalization of this approach. 

Although schemes such as one-stage  \GCSI{k}{\corrector} or  two-stage  \GCSI{k}{\predictor,\corrector} are more costly per time step than GCE, by adding more implicit stages to the scheme we stabilize the method, therefore allowing for significantly larger time steps as compared to the GCE scheme. This improves the overall performance of the solver --- see Section~\ref{sec:numerical-results}, where we discuss the results of numerical experiments.

\begin{algorithm}[H]
\vspace{3ex}

\DontPrintSemicolon
\KwData{$\ur^{n-i},\vr^{n-i},\wr^{n-i}, \quad i = 1,\dots,k$}
\KwResult{$\ur^{n},\vr^{n},\wr^{n}$}

\Begin{
$\wre :=\wr^{n-1}$ \;
 \vspace*{0.5\baselineskip}
\SetKw{kwIn}{in}
\For{$\mode$ \kwIn {$(\mode_1,\ldots,\mode_\nofmodes)$ } } {
$\ure  :=  \gamma\timestep \wre-\sum\limits _{i=1}^{k}\alpha_{i}\ur^{n-k}$ \Comment*{Explicit step}
$\mapae =\Id+ \ure,\qquad \Oe=\mapae(\Or)$ \Comment*{New geometry}
 \vspace*{1\baselineskip}
Find  $\vci \in H^{1}(\Oe)$ and $\implicit{\pressure}\in L^2(\Oe):$ \Comment*{$M$-mode step}
\begin{equation}
\label{eq:ggce:implicit}
\begin{cases}
\begin{array}{lll}
{a}_{\mode}(\vci,\testv)+{b}(\testv,\implicit{\pressure}) & =g_{\vc\mode}(\testv)& \forall\testv\in H_{D}^{1}(\Oe),\\
{b}(\vci,q) & = g_p(q)\,& \forall q\in L^{2}(\Oe).
\end{array}
\end{cases}
\end{equation}\;
$\wri :=\Ext(\vri)$ \Comment*{Extension}
 \vspace*{0.25\baselineskip}
$\wre :=\wri$
 \vspace*{0.25\baselineskip}
 }
 
 \vspace*{0.25\baselineskip}
 
 $\vr^n:=\vri ,\qquad  \wr^{n} :=\wri$\;
 
 $\ur^n:=\gamma\timestep \wr^{n}-\sum\limits _{i=1}^{k}\alpha_{i}\ur^{n-i}$ \Comment*{Displacement recovery}
 }

\vspace{2ex}

\caption{$k$-th order, $\nofmodes$--stage \GCSI{k}{\mode_1,\ldots,\mode_\nofmodes} scheme. The bilinear forms $a_\mode$, $b$ and functionals $g_{\vc\mode}$, and $g_\pressure$ are defined in \eqref{eq:form_a_j},\eqref{eq:form_b} \eqref{eq:form_g_j}, and \eqref{eq:form_g_p}, respectively.
\label{alg:1}
\vspace{3ex}
}
\end{algorithm}

The bilinear forms $a_\mode(\cdot, \cdot)$, $b(\cdot, \cdot)$ and functionals  $g_\mode(\cdot)$ which appear in \eqref{eq:ggce:implicit} are as follows:
\begin{align}
  \label{eq:form_a_j}
  a_{\mode}(\vci,\testv):=&\left(\rho\left(\aleDtau{\vci}+(\grad\vc^{\star}_\mode)\vc^{\circ}_\mode\right),\;\testv\right)_{\Omega^{\#}}+ \left(\eta\epsilon(\vci),\epsilon(\testv)\right)_{\Omega^{\#}},\\
  \label{eq:form_b}
  b(\vci, q):=& \left(\nabla\cdot \vci,  q \right)_{\Omega^{\#}},\\
  \label{eq:form_g_j}
  g_{\vc\mode}(\testv):= &\left(\body\cdot\testv \right)_{\Omega^{\#}} +\left( \traction^{*}\cdot\testv \right)_{\Gamma_{N}^{\#}} \nonumber\\
 & -\mu_{\solid}\left(2\symgrad\left(\sum\limits _{i=1}^{k}\alpha_{i}\uc^{n,i}\right) -(\grad \ucse)^T\grad \ucse,\testv\right)_{\Os^{\#}},
\end{align}
where $\vc^{\circ}_\mode$, $\vc^{\star}_\mode$ are prescribed below. Let us note in the passing that  $\aleDtau{\vci}$ involves velocities $\vci, \vc^{n-1},..., \vc^{n-k}$.

\begin{itemize} 
\item If $\mode=\predictor$, we define $\vc^{\circ}_\mode$ using explicit extrapolation, as in~\cite{van1986second}: depending on the order of the scheme $k$,  
\begin{equation}
\begin{split}\vc^{\circ}_\mode & :=\vc^{n,1}-\wce
\qquad\text{for }k=1,\\
\vc^{\circ}_\mode & :=2(\vc^{n,1}-\wce)-(\vc^{n,2}-\wc^{n,1})\qquad\text{for }k=2.
\end{split}
\label{eq:v_circ_extrapolation}
\end{equation}
For the convective velocity $\vc^{\star}_\mode$ we use an explicit velocity
\begin{equation}
\begin{split}\vc^{\star}_\mode & :=\vc^{n,1}\qquad\text{for }k=1,\\
\vc^{\star}_\mode & :=2\vc^{n,1}-\vc^{n,2}\qquad\text{for }k=2,
\end{split}
\end{equation}
which makes the bilinear part of the form $a_{\mode}(\cdot,\cdot)$
symmetric.
\item For $\mode=\corrector$, we make the step a bit more implicit and set
\begin{equation}
\vc^{\circ}_\mode:=\vc^{\#}-\wce,
\end{equation}
and turn to the semi-implicit advection,
\begin{equation}
\label{eq:v_star}
\vc^{\star}_\mode:=\vci,
\end{equation}
to ensure a better stability of the scheme. For further discussion we refer to \cite{dong2010unconditionally} or \cite{turek1996comparative}. 
Since the velocity is updated after the corrector step, we use the (cheaper) explicit advection in the predictor. This provides a good balance between accuracy and efficiency.
\end{itemize}

The functional $g_p(\cdot)$ will be introduced in the following section, see \eqref{eq:form_g_p}.

\subsubsection{Volume-preserving correction}

\label{sec:vol-stab}
In our case, $\rhos=\text{const}$ and mass conservation is equivalent to incompressibility, i.e.,
\begin{equation}
\div   (\vcs)= -\frac{1}{\rhos}\frac{\partial\rhos}{\partial t} = 0. \label{eq:contraint_solid}
\end{equation}
If, for whatever reason, the solid density is perturbed (e.g., by a numerical approximation of the zero-divergence constraint), such perturbation is accumulated and maintained through time evolution, resulting in a solid volume that may change as a result of these errors. 

We provide a dynamic and strongly consistent correction to the volume of the solid by introducing an additional term in \eqref{eq:contraint_solid} with the aim of restoring $\Jr=1$ whenever the scheme moves away from it:
\begin{equation}
\frac{\partial\rhos}{\partial t}= J\frac{\rhos}{\eta_{V}}\left(\frac{\rhos}{\rho_{s\,0}}-1\right) = J\frac{\rhos}{\eta_{V}}(J-1)\label{eq:volume-new},
\end{equation}
so that the solution would approach the density $\rho_{s\,0}$ regardless
of its starting point. The damping parameter $\eta_{V}$ can be interpreted as a volumetric viscosity which controls the dynamic rate of density correction and leads to a modified volume constraint for the solid with volume-preserving correction,
\begin{equation}
\div  (\vcs)=-\frac{J}{\eta_{V}}(J-1).
\end{equation}
Accordingly, we reformulate the weak form \eqref{eq:weak-form-fsi-ale}${}_2$ of the volume constraint in the solid by adding the respective right-hand side, expressed in the reference configuration as
\begin{equation}
\int\limits _{\Os}q\left(\div   v\right)\df \xc=-\int\limits _{\Ors}\frac{1}{\eta_{V}}(\text{det}(\Fr)-1)\hat{q}\,\text{d}\xr\quad\forall q\in L^{2}(\Omega). \label{eq:solid-damped-div}
\end{equation}
In the actual time-discrete scheme outlined in Algorithm~\ref{alg:1}, the volume-preserving correction is treated in an explicit manner, so that the right-hand side in the weak form \eqref{eq:ggce:implicit} is defined as
\begin{equation}
  \label{eq:form_g_p}
  g_p(q) := \left(-\frac{1}{\eta_{V}}(\Jr^{\#}-1), \hat{q}\right)_{\Ors},
\end{equation}
where $\Jr^{\#}=\det(\Fr^{\#})$, and $\Fr^{\#}=\Id+\Grad\urs^{\#}$ is the most recent algorithmic approximation of the deformation gradient $\Fr$.

\subsection{Spatial discretization\label{sec:Spatial-discretization}}

Let us now introduce the fully discrete approximation of the fluid-structure
interaction problem. We consider triangulation $\Th$ of domain $\Or$
with characteristic element size of $h$. In our case, triangulation
$\Th$ consists of quadrilateral (2D) or hexahedral (3D) elements.
We consider a matching grid, i.e.\ assume
that the initial fluid-solid interface does not intersect with any
element. With sets of polynomials $\mathcal{P}_{p}(K)$ of order
$p$ on each element $K$, we define the finite element spaces
on triangulation $\Th$ 
\begin{equation}
  \label{eq:spaces}
\begin{split}\hat{\mathbb{V}}_{h} & =\{v\in H^{1}(\Or)\,:\,v|_{K}\in\mathcal{P}_{p_{1}}(K)\quad\forall K\in\Th\}^{d},\\
\hat{\mathbb{Q}}_{h} & =\{q\in H^{1}(\Or)\,:\,q|_{K}\in\mathcal{P}_{p_{2}}(K)\quad\forall K\in\Th\}.
\end{split}
\end{equation}
where $p_{1}$ and $p_{2}$ are the orders of finite elements for
the velocity and pressure, respectively. We first discretize the solid
displacement $\urs\in\hat{\mathbb{V}}_{h}$ and the pseudo-displacement
$\ur\in\hat{\mathbb{V}}_{h}$ that defines discrete mapping
$\mapa^{\#}$. We then define the triangulation of $\Omega^{\#}$ as a
transformed triangulation $\Th$ by mapping $\mapa^{\#}$.
Note that triangulation $\Th$ is a matching triangulation
of $\Omega^{\#}$ i.e. the solid-fluid interface does not intersect
with any element. The finite element spaces on domain $\Omega^{\#}$
are defined as: 
\begin{equation}
\begin{split}\mathbb{V}_{h} & =\{v\circ(\mapa^{\#})^{-1}\;:\;v\in\hat{\mathbb{V}}_{h}\}^{d},\\
\mathbb{Q}_{h} & =\{q\circ(\mapa^{\#})^{-1}\;:\;q\in\hat{\mathbb{Q}}_{h}\}.
\end{split}
\end{equation}
We assume that $\mathbb{V}_{h}$ and $\mathbb{Q}_{h}$ satisfy the
Ladyzhenskaya\textendash Babuska\textendash Brezzi condition~\cite{boffi2013mixed}.

\paragraph{Streamline stabilization\label{subsec:Streamline-upwind-stabilization}}

Problem \eqref{eq:form_a_j} is of convection-diffusion type, so 
for  flows with high Reynolds numbers the convection becomes dominating. Since in such cases a straightforward finite element
discretization typically results in oscillatory solutions, some additional stabilization
to the original form $a_j(\cdot,\cdot)$ is necessary; in our scheme, we simply used
\begin{equation}
\label{eq:a_stab}
a^{\text{stab}}_{j}(\vc,\testv)=a_{j}(\vc,\testv)+\sum_{K\in \Th}\int_{K}r \, (v^{\circ}\cdot\grad v) \, (v^{\circ}\cdot\grad\testv)\, \df\xc;
\end{equation}
another possibility would be to use, e.g., the streamline-upwind Petrov--Galerkin (SUPG) scheme~\cite{brooks1982streamline}. The stabilization parameter inside cell $K$ is \cite{brooks1982streamline,john2006discontinuity}
$$
\mathit{\mathit{r}=}\frac{h_K}{2\|v^{\circ}\|p_{1}}\frac{\text{coth}(\Pe_K)-1}{\Pe_K},
$$
where $h_K$ is the diameter of cell $K\in \Th$ and the Peclet number
$\Pe_K$ computed with respect to the cell size is 
$$
\Pe_K=\|v^{\circ}\|\frac{h_K}{2\mu\,p_{1}}.\label{eq:CellPecelet}
$$
In cells where $\Pe_K<1$ we set $r=0$ to avoid problems
with floating-point arithmetics. Note that for a sufficiently
fine mesh $\Pe_K<1$ and thus the stabilization term vanishes,
therefore it does not affect the solution. However, we intend to use a multigrid preconditioner to solve the linear system, thus we also need non-oscillatory solution regardless how large the element size is.

\section{A multilevel, matrix-free preconditioner for the linear system}
\label{sec:preconditioner}

On every time step, there are two computationally intensive parts inside the ``\textbf{for}'' loop in Algorithm~\ref{alg:1}: the computation of the extension of the velocity,  $\Ext(\vri)$  on the fluid domain, and the solution of the discretized system \eqref{eq:ggce:implicit} in order to determine the velocity and the pressure on the entire domain. To perform the former, we use a standard multigrid preconditioned CG, which in our experiments worked just fine. The latter system is much more challenging and therefore we will focus on this subproblem in the present section. 

The linear system \eqref{eq:ggce:implicit} has a block structure of a generalized Oseen-type saddle point problem
\begin{equation}
\label{eq:stokes-discrete-block}
\begin{bmatrix}A & B^{T}\\
B & 0
\end{bmatrix} 
\begin{bmatrix}
\vc \\ p
\end{bmatrix}
= 
\begin{bmatrix}
g_\vc \\ g_p
\end{bmatrix},
\end{equation}
where the square matrix $A$ corresponds to a discretized
convection--diffusion--reaction operator \eqref{eq:a_stab} with the viscosity
coefficient which is discontinuous across the fluid-solid interface, while $B$
corresponds to the discrete divergence operator. Since the number of unknowns in \eqref{eq:stokes-discrete-block} is very large, cf. Table\ref{tab:dofs}, direct solution of this system is infeasible. On the other hand, the system is also challenging for iterative solvers: it is
ill-conditioned due to both the fine mesh size and the high-contrast in the
effective viscosities between the fluid and solid domains ($\eta_{\fluid}$ vs.
$\gamma\tau\mu_{\solid}$), so it requires an efficient preconditioner. Since
the fluid and the solid typically have very different properties, the viscosity
contrast plays a substantial role even when the time step $\timestep$ is
relatively small. %

Many existing preconditioners for FSI problems in the monolithic formulation, e.g. \cite{xu2015well,yang2015modeling}, exploit the block structure of the problem. Others may be based on the multigrid, cf.\ \cite{hron2006monolithic,gee2011truly,richter2015monolithic} or domain decomposition method \cite{crosetto2011parallel,wu2014fully}. See \cite{langer2016recent} for a broad survey of recent developments in this field.

We base our implementation on the \texttt{deal.II}
library~\cite{dealii2019design}, taking advantage of the specific structure of
the linear problem, and leveraging the properties of modern computer hardware,
such as the availability of vectorized  SIMD instructions and parallelism ---
by choosing to implement the preconditioner (and the solver) using the matrix-free approach
\cite{Kronbichler2012}, which is very well supported in the \texttt{deal.II}
library. The CPU cache efficiency is improved because the data is accessed in a more localized manner, reducing the number of cache misses and increasing the overall performance of the solver~\cite{KronbichlerEtAl2022}. Furthermore, the solver's memory footprint is kept low, which means that it can deal with larger problems or function on computers with limited memory resources.

To this end, we adapt the multilevel preconditioner developed for a
generalized stationary Stokes problem with discontinuous viscosity coefficient,
proposed and analyzed in \cite{wichrowski2022}. This preconditioner is not only
robust with respect to the mesh size and the viscosity contrast, but also
supports the matrix-free paradigm by design. Theoretical foundations for a very similar preconditioner have recently been laid in \cite{jodlbauer2022matrixfree}.
While matrix-free preconditioners have successfully been applied in  various
contexts, e.g. phase-field fracture problems \cite{jodlbauer2020fracture} or
mantle convection simulations \cite{kronbichler2012high}, its application
within FSI solver frameworks is much less common, and state-of-the-art FSI
solvers are often matrix-based~\cite{africa2023lifexcfd}. See also
\cite{davydov2020matrix} for a recent application of matrix-free methods in the
case of finite-strain solids.

If the convection is treated explicitly, \eqref{eq:stokes-discrete-block} already defines a generalized Stokes system, so the preconditioner of \cite{wichrowski2022} can be applied directly. In the case of implicit convection, however, the matrix $A$ becomes nonsymmetric, so the preconditioner needs some adjustments. For the convenience of the reader, below we briefly recall the main ingredients of the preconditioner, already adapted to the nonsymmetric case. 

Since the preconditioner is of multilevel type, we will assume that the underlying grid $\Th$ is a result of a $J$--level uniform refinement of some coarse grid $\Th_{0}$ aligned with the fluid-solid interface, resulting in a family of nested conforming triangulations of $\Or$:
\begin{equation}
\mathcal{T}_{0}\subset\mathcal{T}_{1}\subset...\subset\mathcal{T}_{J}=\Th.
\end{equation}
These, in turn, generate a family of discrete problems 
defined on mesh $\mathcal{T}_{j}$ , with corresponding block matrices 
\begin{equation}
\mathcal{M}_{j}=\begin{bmatrix}A_{j} & B_{j}^{T}\\
B_{j} & 0
\end{bmatrix}.
\end{equation}

The preconditioner is formulated as $n$ iterations of the V-cycle multigrid for $\mathcal{M}_{J}$, see \cite{zulehner2000class}, with $m$ smoothing steps which use a customized block smoother on the $j$-th level
\begin{equation}
\K_{j}=\left[\begin{array}{cc}
\hat{A}_{j} & B_{j}^{T}\\
B_{j} & B_{j}\hat{A}_{j}^{-1}B_{j}^{T}-\hat{S}_{j}
\end{array}\right]^{-1},
\label{eq:smoother}
\end{equation}
see \cite{wichrowski2022} for details. In order to apply $\K_j$ to a vector, two solves with $\hat{A}_{j}$ and one with $\hat{S}_{j}$ are required. Both are implemented as matrix-free operators as follows:
$$
\hat{A}^{-1}=\Cheb(A,\diag(A),k_{A}), \qquad \hat{S}^{-1}=\begin{cases}
\ChebMINRES(n_{S},k_{S}) \qquad \text{ when }A=A^T,\\
\ChebBiCGStab(n_{S},k_{S}) \qquad \text{ otherwise}.
\end{cases}
$$
Above, $\Cheb(M,D,k)$ denotes  $k$ iterations of the Chebyshev smoother \cite{ADAMS2003593} preconditioned with $D$. By $\ChebMINRES(n_{S},k_{S})$ we denote the result of $n_{S}$ iterations of the MINRES method applied to the system with matrix $S = B\hat{A}^{-1}B^{T}$, preconditioned with  $\Cheb(S,\diag\left(B(\diag{A})^{-1}B^{T}\right),k_{S})$. The only difference between $\ChebMINRES$ and $\ChebBiCGStab$ is that in the MINRES iteration in the former is replaced with the BiCGStab method in the latter. All these iterative methods take  zero vector as the initial guess.

\section{Numerical results}
\label{sec:numerical-results}

In this section, we present the results of numerical experiments conducted to evaluate the performance and efficiency of the method on a set of benchmark problems. %
As mentioned in Section~\ref{sec:preconditioner}, the experimental framework was implemented usinfg the \texttt{deal.II} library \cite{dealii2019design} and its matrix-free framework \cite{Kronbichler2012}. The stable Taylor--Hood finite element pair, which corresponds to choosing $p_{1}=2$ and $p_{2}=1$ in~(\ref{eq:spaces}), was consistently used in all tests. 

\subparagraph{Test problems}

\begin{figure}
  \begin{center}
    \includegraphics{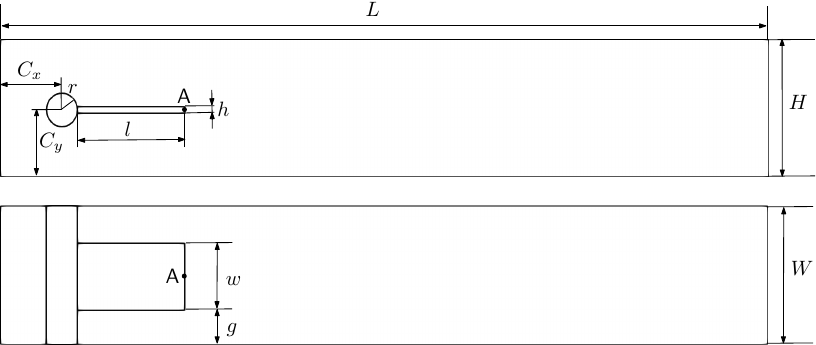}
  \end{center}
  \caption{Geometry for test problems: $x{-}y$ view (top), used for experiments in two dimensions (FSI2i, FSI3i), and as base for the extrusion in three dimensions (FSI3Di), see the $x{-}z$ view (bottom). The plotting point $\mathsf{A}$ is shown in both views.}
  \label{fig:Turek-bench-geometry}
\end{figure}

We evaluate the solver on test problems very closely resembling the well-known Turek--Hron two-dimensional benchmarks FSI2 and FSI3 \cite{turek2006proposal}. These benchmarks consist of an elastic beam attached to a cylindrical obstacle inside a channel, which interacts with the flow. 
To test the solver in 3D we use the benchmark problem FSI3D from  \cite{failer2020parallel}, which consists of an elastic plate attached to a cylinder (depicted in Figure~\ref{fig:Turek-bench-geometry}). This geometry is  obtained by extruding the geometry of 2D tests into the third dimension; the detailed geometrical settings are presented in Figure~\ref{fig:Turek-bench-geometry} together with dimensional parameters in Table~\ref{tab:turek_geometry}.  

Our experimental framework differs from the original benchmarks mentioned above in two ways.
The most important departure is due to the choice of the model of the solid (see Section~\ref{sec:MR-solid}) which, in contrast to Turek and Hron setting, is assumed incompressible. To mark that our tests problems involve incompressible solid, we will refer to them as FSI2i, FSI3i, and FSI3Di, respectively. 
For this reason, our goal will not be to replicate the exact results of the original benchmarks, but rather to demonstrate the performance of our method in simulation of incompressible fluid--structure interaction problems. (To the best of our knowledge, there are no available results for Turek--Hron benchmark problems using an incompressible solid.)

\begin{table}
\caption{Geometry of test problems --- list of dimensional parameters (all expressed in meters).\label{tab:turek_geometry}}
\centering{}%
\begin{tabular}{lccccccccc}
\hline 
$L$ & $H$ &$W$ & $l$ & $h$& $w$ & $C_{x}$ & $C_{y}$ & $r$ &$g$ \tabularnewline
\hline 
2.5 & 0.41 & 0.41 & 0.35 & 0.02& 0.2  & 0.2 & 0.2 & 0.05& 0.1\tabularnewline
\hline 
\end{tabular}
\end{table}

 Secondly, in contrast to the approach used in \cite{turek2006proposal}, where the flow gradually accelerates, our simulation starts with an instantaneous acceleration of the flow at the initial time. This modification allows us to evaluate the solver's capabilities in handling a fully developed flow in a nearly undeformed configuration, and provides a fixed geometry for solver performance comparison --- particularly when considering various time step sizes. However, it induces artificial pressure jumps during the first few time steps. From our observations, it can be seen that these pressure oscillations are damped within a few subsequent initial time steps and have a negligible impact on the final results.

Apart from that, all other settings are identical as in the benchmarks mentioned above. 

In order to improve mass conservation and incompressibility of the solid, unless stated otherwise, we will apply the volume-preserving correction, cf. Section~\ref{sec:Time-derivative-approximation}, with default parameter $\eta_{V}=0.1$.

\subparagraph{Mesh and deformation handling}

\begin{figure}
  \begin{center}
    \includegraphics[width=\textwidth]{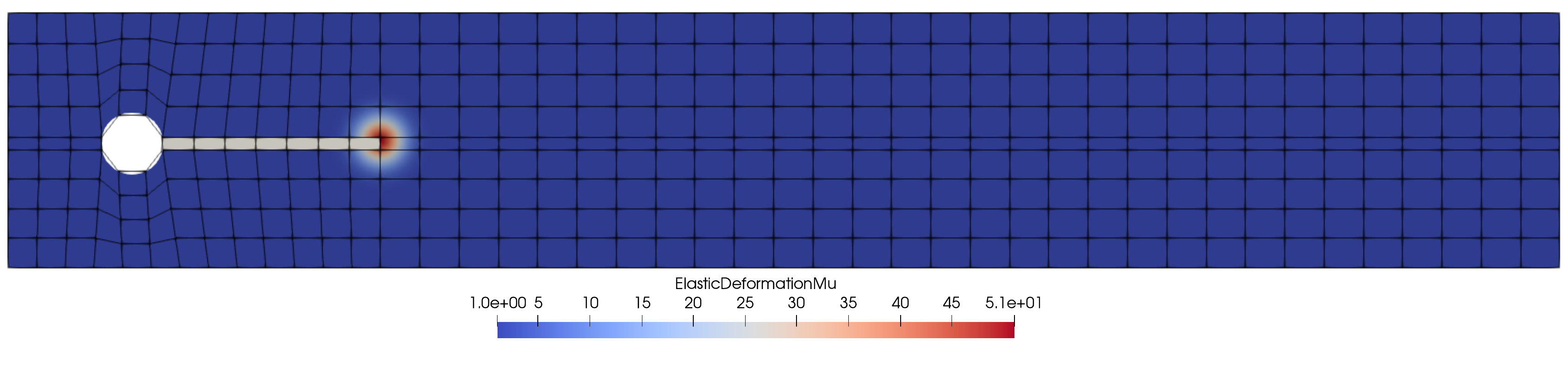}
  \end{center}
  \caption{Coarse grid and distribution of the mesh stiffness parameter, see \eqref{eq:mu_distribution}. Solid marked in grey.}
  \label{fig:coarse_grid}
\end{figure}

To discretize the reference spatial domain $\Or$, we start with a coarse grid $\Th_0$ as illustrated in Figure~\ref{fig:coarse_grid}, which then undergoes $J$ rounds of uniform refinement, including adjustments to accommodate the curved boundary of the cylinder. This results in a family of refined grids $\Th_0,\ldots,\Th_J$. The numbers of degrees of freedom corresponding to selected unknowns for specific $J$ used in the experiments are summarized in Table~\ref{tab:dofs}. For brevity, in what follows we denote by $N$ the total number of unknowns in the momentum equation.

\begin{table}
\centering
\caption{Number of degrees of freedom corresponding to grids in 2D or 3D obtained after $J$ refinements of the coarse meshes.  ($N=$ velocity + pressure).
\label{tab:dofs}}
\resizebox{\textwidth}{!}{
\begin{tabular}{lr|cccccccc}
& $J$ & 1 & 2 & 3 & 4 & 5 & 6 & 7 & 8\\
\hline
2D case  & velocity & 12k & 50K & 198k & 790k   & 3.15M    &  12.6M  &   50.4M & 201M\\
& pressure & 1.6k & 6k & 25k & 99k &  395k & 1.58M &  6.30M &  25.2M\\
& $N$ & 14k & 56k & 223k & 889k &  3.55M & 14.2M &  56.7M  & 227M\\
& displacement & 12k & 50K & 198k & 790k   & 3.15M    &  12.6M  &   50.4M & 201M\\
\hline
3D case & velocity & 533k & 4.11M & 32.2M & 255M & -- & -- & -- & --\\
& pressure & 24k & 178k &  1.37M & 10.7M &  -- & -- & -- & -- \\
& $N$  & 557k & 4.27M &  33.6M & 266M &  -- & -- & -- & -- \\
& displacement  & 533k & 4.11M &  32.2M & 255M &  -- & -- & -- & -- \\
\end{tabular}
}
\end{table}

 To handle the mesh deformation, we employ a linear elasticity problem with a variable coefficient distribution to enhance the quality of the resulting computational grid. The coefficient distribution is as follows:
\begin{equation}
\mu_{A}(x,y,z)=1+50\exp(-800((x-A_{x})^{2}+(y-A_{y})^{2})). 
\label{eq:mu_distribution}
\end{equation}

\subparagraph{Linear solver settings}

The solution to~(\ref{eq:ggce:implicit}) was obtained using the FGMRES iterative solver, preconditioned with the method outlined in Section~\ref{sec:preconditioner}, whose parameters are specified in Table~\ref{tab:NS-solver-parameters-1}. The iteration was terminated when the Euclidean norm of the residual dropped below the threshold value $\varepsilon=10^{-6}$ in FSI2i and FSI3i testcases or $\varepsilon=10^{-3}$ in FSI3Di.
The initial guess for the FGMRES iteration was set equal to the most recent values of the velocity and pressure. 

The extension \eqref{eq:extenstion} was computed by the CG iteration preconditioned with a classical matrix-free multigrid method with fourth order Chebyshev smoother based on the diagonal, cf.~\cite{Kronbichler2012}. There, the stopping criterion was always to reduce the Eucledean norm of the residual below~$10^{-6}$.

\begin{table}
\caption{Material parameters.\label{tab:FSI_parameters}}
\centering{}%
\begin{tabular}{lcccc}
\hline 
Parameter & Units & FSI2i & FSI3i &FSI3Di \tabularnewline
\hline 
$\rhos$ & \kilogram\usk\meter\rpcubed %
& $10^{4}$ & $10^{3}$  & $10^{3}$\tabularnewline
$\rhof$ &  \kilogram\usk\meter\rpcubed %
& $10^{3}$ & $10^{3}$& $10^{3}$  \tabularnewline
$\mu_{\solid}$ & \kilogram\usk\reciprocal\meter\usk\second\rpsquared %
& $0.5\times10^{6}$ & $2\times10^{6}$ & $2\times10^{6}$ \tabularnewline
$\eta_{\fluid}$ & \kilogram\usk\reciprocal\meter\usk\reciprocal\second %
& 1 & 1 & 1\tabularnewline
$V_{in}$ & \meter\usk\reciprocal\second %
& 1 & 2  & 1.75  \tabularnewline
$\Rey$ & --- & 100 & 200 & 175\tabularnewline
\hline 
\end{tabular}
\end{table}

\begin{table}
\caption{Default parameters of the preconditioner described in Section~\ref{sec:preconditioner}.\label{tab:NS-solver-parameters-1}}
\centering{}%
\begin{tabular}{llccc}
\hline 
\multicolumn{2}{c}{Parameter} %
&  & $\predictor$ & $\corrector$\tabularnewline
\hline 
Order of the Chebyshev smoother defining $\hat{A}^{-1}$ & $k_{A}$ &  & 4 & 6\tabularnewline
Number $m$ of outer MG smoothing steps & $m$ &  & 2 & 2\tabularnewline
Number of outer MG iterations & $n$ &  & 1 & 1\tabularnewline
Order of the Chebyshev smoother defining $\tilde{S}^{-1}$ & $k_{S}$ &  & 1 & 1\tabularnewline
Number of inner MINRES/BiCGStab iterations & $n_S $ &  & 1 & 1\tabularnewline
\hline 
\end{tabular}
\end{table}

\subsection{Tests in 2D}

In both FSI2i and FSI3i testcases, the cylinder, upper and lower sides of the channel are considered as rigid walls and no-slip boundary conditions are imposed there. At the right end of the channel, a %
do--nothing 
boundary condition is assumed, while on the left side of the channel a parabolic inflow velocity profile is prescribed, cf.~\cite{turek2006proposal}. The flow is defined by the average velocity at the inlet, $V_{in}$, and by the Reynolds number (computed with respect to the obstacle diameter). Here, we consider $V_{in}=1$ or $V_{in}=2$, corresponding to $\Rey=100$ and $\Rey=200$, respectively. The detailed test data is presented in Table~\ref{tab:FSI_parameters}.

We conduct FSI2i and FSI3i tests using the %
\GCSIEI{} scheme on a fixed spatial grid after $J=5$ levels of refinement (see Table~\ref{tab:dofs} for details on the number of degrees of freedom). By varying time step size~$\timestep$ we evaluate the stability and convergence (in time) of our numerical method. %
Additionally, for qualitative comparison, we run the FSI2i benchmark using a one-stage method, \GCSII. Let us note that the standard GCE scheme (which, in our notation, corresponds to the \GCSI{2}{\predictor} method) turned out unstable for time steps larger than $\timestep=0.00025$ (and thus for all time steps considered below), so we excluded this method from the following analysis.

In Tables \ref{tab:FSI2_comparison} and \ref{tab:FSI3_comparison-1} we report the amplitude, average displacement, and period of the last computed oscillation at point $\mathsf{A}$, see Figure~\ref{fig:Turek-bench-geometry}. Apparently, with $\tau$ decreasing, both \GCSIEI{} and \GCSII{} schemes  converge to similar results, with the former settling for larger time step values than the latter. Lower amplitude of oscillations of  $u_{y}(\mathsf{A})$ in the case of \GCSII{} suggests that the single-stage method introduces more artificial viscosity compared to the two-stage scheme.

As expected, the results obtained from FSI2i and FSI3i tests are close to FSI2 and FSI3 from \cite{turek2006proposal}. However, they are not identical, since our solid model uses partially different constitutive laws. 

\begin{table}
\caption{Comparison of displacements and frequency in FSI2i (incompressible material) test problem obtained with \GCSIEI{} (top) and \GCSII{} (bottom) schemes, for various time step sizes $\timestep$. For qualitative comparison, in the last row we also report results of the FSI2 benchmark for a compressible material. Spatial mesh with $J=5$ levels.
\label{tab:FSI2_comparison}}
\centering{}%
\begin{tabular}{llccc}
\hline %
& $\timestep$ & $u_{x}(\mathsf{A})\times10^{-3}$ & $u_{y}(\mathsf{A})\times10^{-3}$ & Frequency  ($u_{y}(\mathsf{A})$)\tabularnewline
\hline
FSI2i, \GCSIEI
& $0.04$  & \multicolumn{3}{c}{(no steady oscillations)} \tabularnewline
& $0.02$  & $-14.83\pm11.65$ & $1.29\pm76.38$ & $1.96$\tabularnewline
& $0.01$   & $-15.09 \pm 13.04$ & $1.25\pm80.26$ & $2.00$\tabularnewline
& $0.005$  & $-14.85\pm12.89$ & $1.23\pm80.77$ & $2.00$\tabularnewline
& $0.0025$  & $-14.72\pm12.55$ & $1.23\pm80.77$ & $2.00$\tabularnewline
\hline %
FSI2i, \GCSII
& $0.04$  & \multicolumn{3}{c}{(no steady oscillations)} \tabularnewline
& $0.02$  & $-13.87\pm11.80$ & $1.38 \pm 76.45$ & $1.92$\tabularnewline
& $0.01$   & $-14.26 \pm 12.65$ & $1.32\pm79.42$ & $1.96$\tabularnewline
& $0.005 $  & $ -14.53 \pm 12.71 $ & $1.30\pm80.51$ & $1.98$\tabularnewline
& $0.0025$  & $-14.52\pm 12.52$ & $1.23\pm80.77$ & $2.00$\tabularnewline
\hline
FSI2 (compressible) \cite{turek2006proposal}
& 0.001 & $-14.58\pm12.44$ & $1.23\pm80.6~$ & $2.0~$ \tabularnewline
\hline 
\end{tabular}

\end{table}

\begin{figure}
    \includegraphics{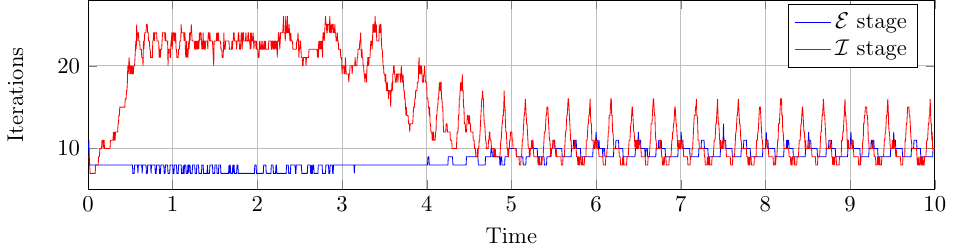}

        \caption{Number of FGMRES iterations per time step in 2D case (FSI2i), as a function of time. Results for $J=5$ and $\timestep = 0.005$. \label{fig:FSI2iter} }
  \end{figure}

The number of FGMRES iterations required to reduce the residual below prescribed threshold remains, as seen in Figure~\ref{fig:FSI2iter}, roughly constant over time, with lower convergence rate observed during the initial stages of the simulation, %
which then improves in the steady oscillation phase. During this second phase, the solver performs no more than 15 iterations in $\corrector$ stage of the \GCSIEI{} scheme (and significantly less in $\predictor$ mode), confirming high efficiency of the preconditioner described in Section~\ref{sec:preconditioner} even when the linear problem is nonsymmetric (see also Table~\ref{tab:FSI3D-iter-J}). %

\begin{table}
\caption{Comparison of displacements and frequency in the FSI3i (incompressible) testcase with $J=4$. For qualitative comparison, in the last row we also report results of the FSI3 benchmark for a compressible material. \label{tab:FSI3_comparison-1}}
\centering{}%
\begin{tabular}{llccc}
\hline 
& $\timestep$ & $u_{x}(\mathsf{A})\times10^{-3}$ & $u_{y}(\mathsf{A})\times10^{-3}$ & Frequency ($u_{y}(\mathsf{A})$)\tabularnewline
\hline
FSI3i, \GCSIEI
& $0.02$  & \multicolumn{3}{c}{(no steady oscillations)} \tabularnewline
& $0.01$ & $-2.85\pm2.06$ & $1.80\pm32.19$ & $5.38$\tabularnewline
& $0.005$ & $-2.79\pm2.46$ & $1.47\pm34.38$ & $5.47$\tabularnewline
&  $0.0025$ & $-2.88\pm2.69$ & $1.45\pm34.56$ & $5.51$\tabularnewline
\hline
FSI3 (compressible) \cite{turek2006proposal}
& 0.0005 &  $-2.69\pm2.53$ & $1.48\pm34.38$ & $5.3~$\tabularnewline
\hline 
\end{tabular}
\end{table}

\pgfplotstableread[comment chars={"},col sep=comma]{./results/FSI2/GCE0/J51e-2-GCE0.csv}\FSInoCorrector

\pgfplotstableread[comment chars={"},col sep=comma]{./results/FSI2/GCE2/Full_solve/J51e-2-full.csv}\FSIwithCorrector

\subsection{Tests in 3D }

In the FSI3Di test case, a non-slip condition is imposed  on the long outer sides of the channel. On the right side of the domain a do-nothing outflow condition is precsribed, while the left side is subject to the prescribed inflow velocity profile Dirichlet boundary condition:
\begin{equation}
v_{in}(x,y,z) = \frac{36 V_{in} } { H^2W^2} y (H - y) z (W - z),
\end{equation}
which results in $\Rey=175$ for $V_{in}=1.75$  (calculated with respect to the cylinder diameter). For the handling of grid deformations, we use the linear elasticity problem with the same coefficient distribution as in 2D case, cf.~(\ref{eq:mu_distribution}).

To generate the computational grids, we first extrude the 2D grid used in our earlier experiments along the third dimension, resulting in 6 elements in that direction. Then, the grid is uniformly refined $J$ times to obtain a multilevel structure. We conduct the experiments dealing with 4 different levels of refinements, $J=1, \ldots,4$ (see Table~\ref{tab:dofs} for details on the number of degrees of freedom) and several time step sizes.  

To assess the quality of the solution, we again record the oscillation
frequency, the average, and the amplitude of the displacement of point
$\mathsf{A}$, see Table~\ref{tab:FSI3D_data}. The results indicate that, as the
time step size decreases, the oscillation frequency, average and amplitude of
displacement of point $\mathsf{A}$, each settle at a certain value, supporting
the expectation that the scheme is convergent as~$\tau\to 0$. 

As in the 2D case, the results of the FSI3Di test are again slightly different from FSI3D~\cite{failer2020parallel}, %
due to different material laws for the solid. %

\begin{figure}
    \includegraphics{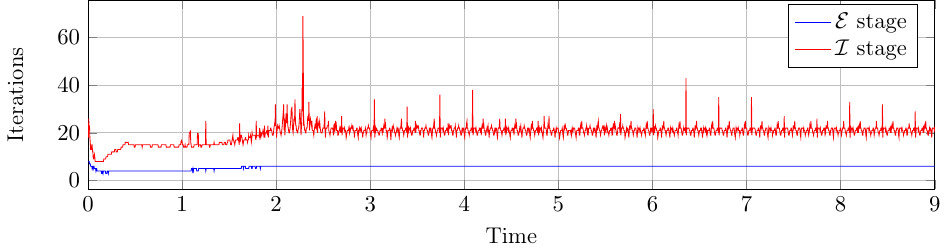}
    \caption{Number of FGMRES iterations per time step in 3D case (FSI3Di), as a function of time. \label{fig:FSI3Diter}}
  \end{figure}

From Figure~\ref{fig:FSI3Diter} it follows that the number of FGMRES iterations required to converge is almost always low, or very low, respectively, during the $\corrector$ or $\predictor$ stage of the \GCSIEI{} scheme. It also remains roughly constant over time. This shows that the preconditioner works efficiently for the linear systems of equations being solved in each stage of the \GCSIEI{} scheme in~3D as well; see also Table~\ref{tab:FSI2D-iter-J}.

\begin{figure}
\begin{centering}
\subfloat[Geometry\label{fig:FSI3D-geometry}]{\centering{}\includegraphics[width=0.48\columnwidth]{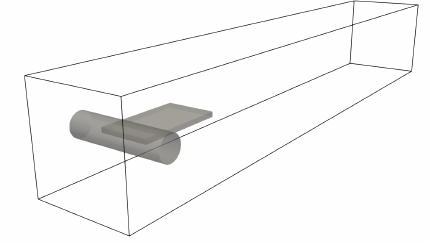}}
\subfloat[Streamlines\label{fig:FSI3D-streamlines}]{\centering{}\includegraphics[width=0.48\columnwidth]{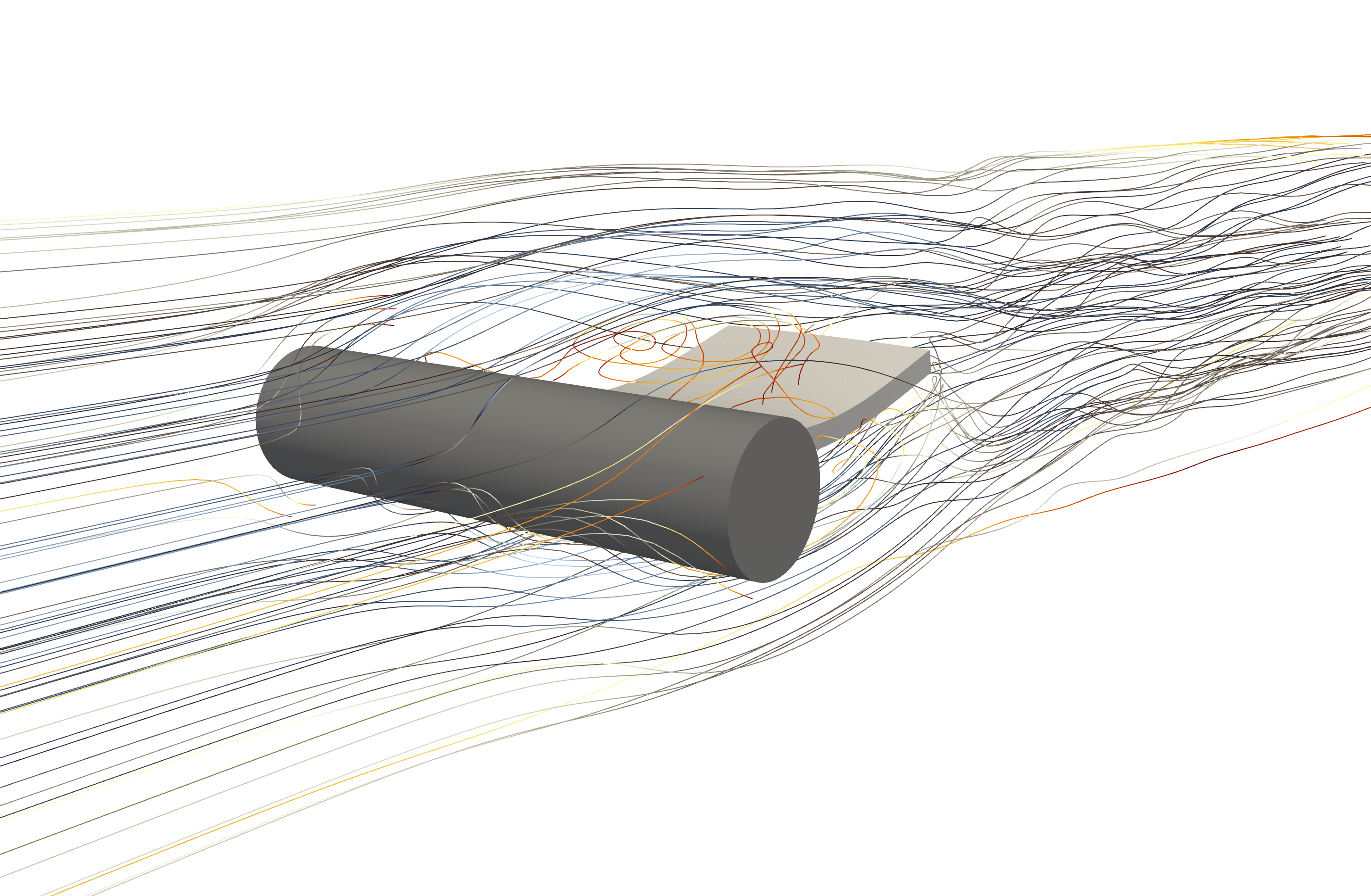}}

\end{centering}
\centering{}\caption{(a) Visualization of the geometry of the FSI3Di testcase. (b) Plot of the streamlines and of the deformed plate.
\label{fig:FSI3D-turek}}
\end{figure}

\begin{table}
\caption{Comparison of displacements and frequency in FSI3Di test with $J=3$. For qualitative comparison, in the last row we also report results of the FSI3D benchmark for a compressible material.  \label{tab:FSI3D_data}
}
\centering{}%
\begin{tabular}{ccccc}
\hline 
& $\timestep $& $u_{x}(\mathsf{A})\times10^{-3}$ & $u_{y}(\mathsf{A})\times10^{-3}$ & Frequency  $u_{y}(\mathsf{A} )$)\tabularnewline
\hline
 FSI3i, \GCSIEI
 & $0.008$  & $ -2.39 \pm 2.23 $ & $ 2.51 \pm  28.47$ & $5.63 $\tabularnewline
 & $0.004$  & $ -2.39 \pm 2.44  $ & $ 2.45 \pm 29.20 $ & 5.73\tabularnewline
 & $0.002$  & $- 2.39  \pm    2.49 $ & $ 2.43  \pm  29.24 $ & $ 5.77 $\tabularnewline
\hline
{l}{FSI3D (compressible)} \cite{failer2020parallel} & $-2.143 \pm 2.383 $ & $ 2.699 \pm 25.594 $ & $5.60$\tabularnewline
\hline 
\end{tabular}
\end{table}

\pgfplotstableread[comment chars={"},col sep=comma]{./Iter/SLURMS/csvs/FSI2/FSI2-5e3.csv}\JiterFull
  
\pgfplotstableread[comment chars={"},col sep=comma]{./FSI3D_ref/30M-4e-3.csv}\threeDiter

\subsection{Influence of the volume-preserving correction on the quality of the solution}%

\begin{figure}
  \includegraphics{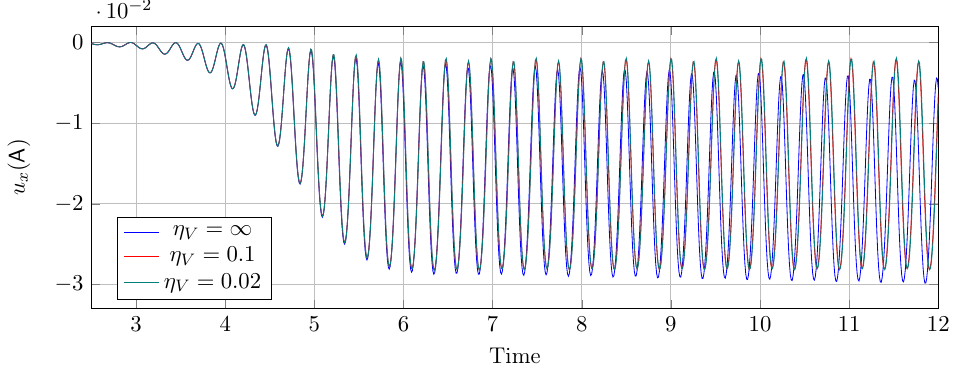}
    
 \caption{Influence of volume-preserving correction on the displacement of the point $\mathsf{A}$ versus time, for varying damping parameter $\eta_V$ of the volume-preserving correction. (By convention, $\eta_V = \infty$ corresponds to no correction.) FSI2i test case, grid with $J=5$ levels, \GCSIEI{} scheme with timestep $\timestep=0.01$. 
\label{fig:FSI-neta-A}}

  \end{figure}

Here we investigate how effective is the volume-preserving correction, introduced in Section~\ref{sec:vol-stab}, in resolving certain issues related to the deformation of the beam. 
In Figure~\ref{fig:FSI-neta-A}, we compare the results obtained for the FSI2i testcase with a relatively large time-step size, $\timestep=0.01$, when the volumetric damping is either weak, strong or absent. It is evident that without the volumetric correction ($\eta_{V}=\infty$) the horizontal displacement $u_{x}$ at point $\mathsf{A}$ is creeping, whereas it remains stable for both $\eta_{V}=0.1$ and $\eta_{V}=0.02$. Moreover, the displacement graphs obtained for the latter values of $\eta_{V}$ practically overlap.

\begin{figure}
  \includegraphics{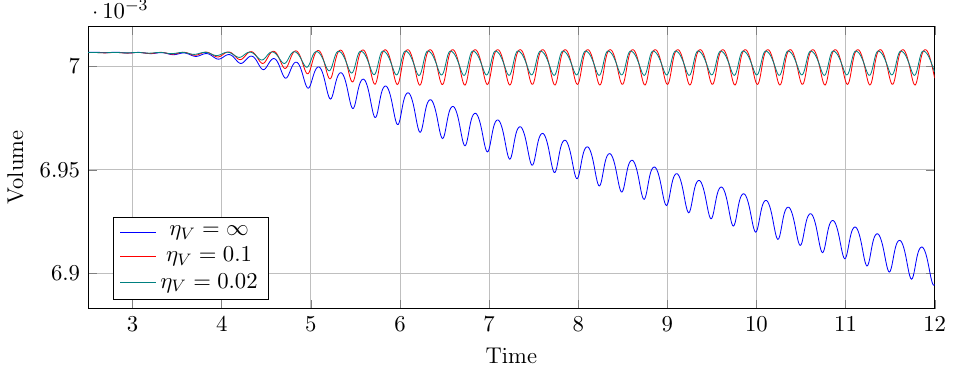}
    
 \caption{Volume of the solid versus time, for varying damping parameter $\eta_V$ of the volume-preserving correction. (By convention, $\eta_V = \infty$ corresponds to no correction.) %
FSI2i test case, grid with $J=5$ levels; \GCSIEI{} scheme with time step $\timestep=0.01$. \label{fig:FSI-volume-damping}}

  \end{figure}

The importance of the volume-preserving correction is even more pronounced when one compares the evolution of the volume of the solid in time, see Figure~\ref{fig:FSI-volume-damping}. Without stabilization, its volume is diminishing (a similar effect also occurs in the case of the first-order time integration scheme \GCSI{1}{\predictor,\corrector} where the beam is gaining volume). The presence of the volume-preserving correction term resolves the problem with the incompressibility condition, thus paving the way to use large time steps  in our scheme. Let us mention that, to some extent, the final result seems quite insensitive to the (small) value of the damping parameter $\eta_{V}$. %

\pgfplotstableread[comment chars={T},col sep=comma]{./results/FSI3/T5e-3/FSI3J5-05-no_vol_damping.csv}\FSInoVolDamping

\pgfplotstableread[comment chars={T},col sep=comma]{./results/FSI2/eta_comparison.csv}\FSIDamping

\pgfplotstableread[comment chars={T},col sep=comma]{./results/FSI3/T5e-3/FSI3J5-05-full.csv}\FSIreference

\pgfplotstableread[comment chars={T},col sep=comma]{./results/FSI2/eta_comparison_volume.csv}\FSIVolume

\subsection{Performance}

\begin{figure} 
  \pgfplotstableread[header=false, comment chars={T},col
sep=comma]{./FSI/performance/FSI2_elaborated.csv}\FSITwoPerformance

\usetikzlibrary{patterns}
  \includegraphics{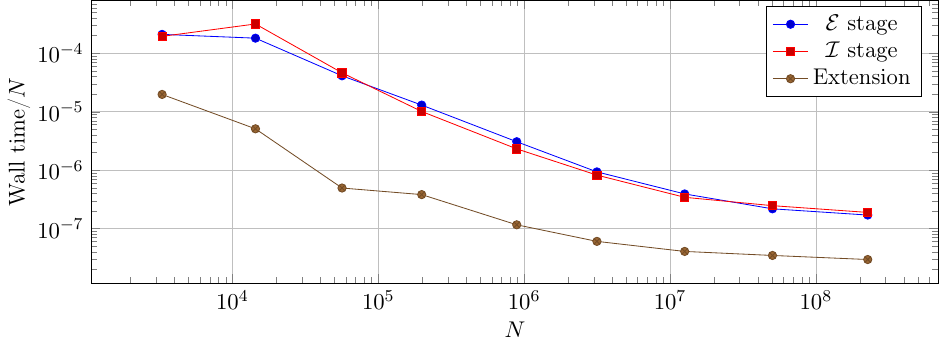}

  \pgfplotstableread[header=false, comment chars={T},col
sep=comma]{./FSI/performance/FSI3D_elaborated.csv}\FSIThreePerformance

\includegraphics{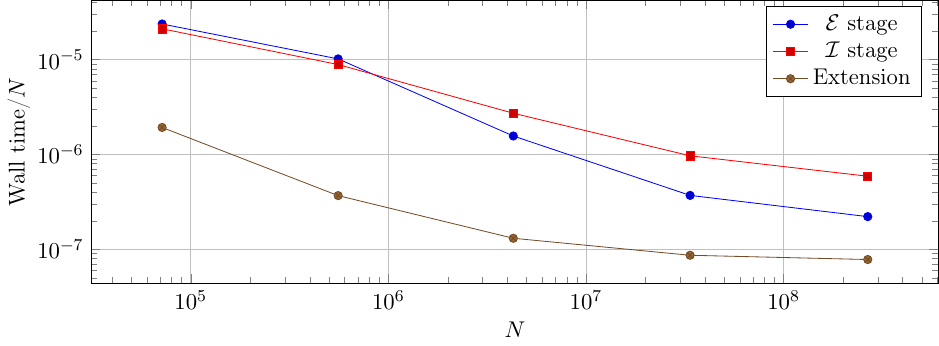}
\caption{Efficiency of the solver, as a function of the total
number $N$ of degrees of freedom in the momentum equation.  FSI2i
problem (top) %
and FSI3Di (bottom). %
\label{fig:FSI-perf}}
\end{figure}

In order to evaluate the performance of our solver,  we plotted the time spent per degree of freedom versus the number of degrees of freedom for both the 2D and 3D cases, as seen in Figure~\ref{fig:FSI-perf}.  The tests were perfomed on 4 nodes, each equipped with two Intel Xeon 8160  @2.1GHz processors. The results demonstrate that our method is highly efficient, with the time spent per degree of freedom decreasing as the number of degrees of freedom increases. This trend can be attributed to the relatively inefficient coarse grid solver becoming a smaller fraction of the overall timing as the problem size increases. We also note that the number of iterations of the FGMRES linear solver required for convergence is roughly independent of the problem size, as demonstrated in Tables~\ref{tab:FSI3D-iter-J}~and~\ref{tab:FSI2D-iter-J}.

\begin{table}
\caption{Number of FGMRES iterations required for convergence ($\varepsilon=10^{-6}$)
at the first time step $\timestep=0.004$ in 2D tests during both stages $\{\predictor,\corrector\}$ of the \GCSIEI{} scheme and for varying number of levels~$J$.\label{tab:FSI3D-iter-J}}
\begin{center}
\begin{tabular}{cccccccc}
\hline
& stage $\downarrow$ $\quad J\rightarrow$& 3 & 4 & 5 & 6 & 7 & 8\tabularnewline
\hline 
\multirow{2}{*}{FSI2i} & $\predictor$ & 8 & 9 & 9 & 10 & 10 & 10\tabularnewline
 & $\corrector$ & 15 & 13 & 11 & 10 & 9 & 9 \tabularnewline
\hline 
\multirow{2}{*}{FSI3i} & $\predictor$ & 11 & 11 & 12 & 12 & 12 & 12 \tabularnewline
 & $\corrector$ & 26 & 28 & 23 & 20 & 19 & 18 \tabularnewline
\hline
\end{tabular}
\end{center}
\end{table}  

\begin{table}
\caption{Number of FGMRES iterations required for convergence ($\varepsilon=10^{-3}$)
at the first time step $\timestep=0.005$ in the 3D test during both stages $\{\predictor,\corrector\}$ of the \GCSIEI{} scheme and for varying number of levels~$J$.\label{tab:FSI2D-iter-J}}
\begin{center}
\begin{tabular}{cccccc}
\hline
& stage $\downarrow$ $\quad J\rightarrow$& 1 & 2 & 3 & 4 \tabularnewline
\hline 
\multirow{2}{*}{FSI3Di} & $\predictor$ & 9 & 10 & 8 & 8 \tabularnewline
 & $\corrector$ & 9 & 15 & 18 & 15  \tabularnewline
\hline
\end{tabular}
\end{center}
\end{table}

\FloatBarrier
\section{Conclusions}
\label{sec:conclusions}
In this paper we presented a numerical method for solving fluid--structure interaction (FSI) problems with incompressible, hyperelastic Mooney--Rivlin solid. For a fully-coupled finite element discretization, we have designed a semi-implicit $k$-th order BDF-based time integration scheme, \GCSI{k}{\cdot}, augmented with a correction term aimed at improving the volume preservation. The method uses the arbitrary Lagrangian--Eulerian (ALE) formulation to track the moving parts. 

The scheme, which resembles an $S$-stage predictor--corrector method, 
possesses improved stability properties, as compared to the GCE scheme introduced in \cite{crosetto2011parallel}. This makes it possible to use much larger time steps while still capturing the dynamics of the interaction between the fluid and the solid. 
Although each time step is roughly $S$ times more expensive than the corresponding step of the GCE method, numerical experiments demonstrate that already for $S=2$, the \GCSIEI{} scheme is stable enough so that the gains in the efficiency, due to taking larger time steps, easily outweigh the increased cost of a single step.

The splitting between the explicit and implicit parts of the integration step has been designed in such a way that, in the latter, a solution to a generalized nonsymmetric Stokes-like problem needs to be computed. %
This is challenging, because the problem's condition number is adversely
affected not only by the number of the unknowns, but also by a very large
contrast in the coefficients of the underlying Stokes-like PDE. To address this
issue, we adapted a preconditioning method from \cite{wichrowski2022} and
proved its robustness and efficiency in this type of application. While the
number of preconditioned FGMRES iterations fluctuated in time, it was always
bounded by a reasonably moderate constant, regardless of the problem size. 

The method described above has been implemented in a matrix-free fashion using
the \texttt{deal.II} library and performed very well on classical FSI benchmarks
in 2D and~3D, with as many as 250M spatial degrees of freedom, being solved in
parallel on four computing nodes, each equipped with two Intel Xeon 8160
@2.1GHz  CPU processors. 

We believe that this approach has the potential for further improvement -- that
we plan to investigate in future research -- for example through the
introduction of a more efficient coarse solver, which is one of the bottlenecks
in the current implementation, as shown by the efficiency results.

\section{Acknowledgements}
This research has been partially supported by the National Science Centre (NCN) in Poland through Grant No.\ 2015/19/N/ST8/03924
and the EffectFact project (No.\ 101008140) funded within the H2020 Programme, MSC Action RISE-2022.

\bibliographystyle{unsrt}
\addcontentsline{toc}{section}{\refname}
\bibliography{bibliografia}

\begin{thebibliography}{10}

\bibitem{taylor1952fsi}
G.~Taylor.
\newblock Analysis of the swimming of long and narrow animals.
\newblock {\em Proceedings of the Royal Society of London. Series A.
  Mathematical and Physical Sciences}, 214(1117):158--183, 1952.

\bibitem{fernandez2005newton}
M.~{\'A}. Fern{\'a}ndez and M.~Moubachir.
\newblock {A Newton method using exact Jacobians for solving fluid--structure
  coupling}.
\newblock {\em Computers \& Structures}, 83(2):127--142, 2005.

\bibitem{badia2008splitting}
S.~Badia, A.~Quaini, and A.~Quarteroni.
\newblock Splitting methods based on algebraic factorization for
  fluid--structure interaction.
\newblock {\em SIAM Journal on Scientific Computing}, 30(4):1778--1805, 2008.

\bibitem{degroote2009performance}
J.~Degroote, K.~Bathe, and J.~Vierendeels.
\newblock Performance of a new partitioned procedure versus a monolithic
  procedure in fluid--structure interaction.
\newblock {\em Computers \& Structures}, 87(11-12):793--801, 2009.

\bibitem{crosetto2011parallel}
P.~Crosetto, S.~Deparis, G.~Fourestey, and A.~Quarteroni.
\newblock Parallel algorithms for fluid--structure interaction problems in
  haemodynamics.
\newblock {\em {SIAM} Journal on Scientific Computing}, 33(4):1598--1622,
  January 2011.

\bibitem{crosetto2011fluid}
P.~Crosetto, P.~Reymond, S.~Deparis, D.~Kontaxakis, N.~Stergiopulos, and
  A.~Quarteroni.
\newblock Fluid--structure interaction simulation of aortic blood flow.
\newblock {\em Computers \& Fluids}, 43(1):46--57, 2011.

\bibitem{benra2011comparison}
F.-K. Benra, H.~Josef Dohmen, J.~Pei, S.~Schuster, and B.~Wan.
\newblock A comparison of one-way and two-way coupling methods for numerical
  analysis of fluid-structure interactions.
\newblock {\em Journal of applied mathematics}, 2011, 2011.

\bibitem{bucelli2022fsi}
Michele Bucelli, Luca Dede, Alfio Quarteroni, and Christian Vergara.
\newblock Partitioned and {Monolithic} {Algorithms} for the {Numerical}
  {Solution} of {Cardiac} {Fluid}-{Structure} {Interaction}.
\newblock {\em Communications in Computational Physics}, 32(5):1217--1256, June
  2022.

\bibitem{clevenger2019comparison}
T.~C Clevenger and T.~Heister.
\newblock Comparison between algebraic and matrix-free geometric multigrid for
  a stokes problem on adaptive meshes with variable viscosity.
\newblock {\em Numerical Linear Algebra with Applications}, 28(5):e2375, 2021.

\bibitem{failer2020parallel}
L.~Failer and T.~Richter.
\newblock A parallel {N}ewton multigrid framework for monolithic
  fluid--structure interactions.
\newblock {\em Journal of Scientific Computing}, 82(2):1--27, 2020.

\bibitem{bazilevs2013computational}
Y.~Bazilevs, K.~Takizawa, and T.~Tezduyar.
\newblock {\em Computational fluid--structure interaction: methods and
  applications}.
\newblock John Wiley \& Sons, 2013.

\bibitem{richter2017fluid}
T.~Richter.
\newblock {\em Fluid--structure interactions: models, analysis and finite
  elements}, volume 118.
\newblock Springer, 2017.

\bibitem{donea2017arbitrary}
J.~Donea, A.~Huerta, J.-P. Ponthot, and A.~Rodr{\'\i}guez-Ferran.
\newblock Arbitrary {L}agrangian--{E}ulerian {M}ethods.
\newblock {\em Encyclopedia of Computational Mechanics Second Edition}, pages
  1--23, 2017.

\bibitem{richter2010finite}
T.~Richter and T.~Wick.
\newblock {Finite elements for fluid--structure interaction in ALE and fully
  Eulerian coordinates}.
\newblock {\em Computer Methods in Applied Mechanics and Engineering},
  199(41):2633--2642, 2010.

\bibitem{ryzhakov2010monolithic}
PB~Ryzhakov, Riccardo Rossi, SR~Idelsohn, and Eugenio Onate.
\newblock A monolithic lagrangian approach for fluid--structure interaction
  problems.
\newblock {\em Computational mechanics}, 46:883--899, 2010.

\bibitem{dunne2006eulerian}
Th. Dunne.
\newblock {An Eulerian approach to fluid--structure interaction and
  goal-oriented mesh adaptation}.
\newblock {\em International Journal for Numerical Methods in Fluids},
  51(9-10):1017--1039, 2006.

\bibitem{glowinski1994fictitious}
R.~Glowinski, T.-W. Pan, and J.~Periaux.
\newblock A fictitious domain method for external incompressible viscous flow
  modeled by navier-stokes equations.
\newblock {\em Computer methods in applied mechanics and engineering},
  112(1-4):133--148, 1994.

\bibitem{BoffiGastaldiHeltai-2018-a}
D.~Boffi, L.~Gastaldi, and L.~Heltai.
\newblock {\em A Distributed Lagrange Formulation of the Finite Element
  Immersed Boundary Method for Fluids Interacting with Compressible Solids},
  pages 1--21.
\newblock Springer International Publishing, Cham, 2018.

\bibitem{Peskin2002}
Charles~S Peskin.
\newblock {The immersed boundary method}.
\newblock {\em Acta Numerica}, 11(1):479--517, jan 2002.

\bibitem{BoffiGastaldiHeltai-2008-a}
D.~Boffi, L.~Gastaldi, L.~Heltai, and Ch. Peskin.
\newblock On the hyper-elastic formulation of the immersed boundary method.
\newblock {\em Computer Methods in Applied Mechanics and Engineering},
  197(25-28):2210--2231, 2008.

\bibitem{hron2006monolithic}
J.~Hron and S.~Turek.
\newblock {A monolithic FEM/multigrid solver for an ALE formulation of
  fluid--structure interaction with applications in biomechanics}.
\newblock In {\em Fluid--structure interaction}, pages 146--170. Springer,
  2006.

\bibitem{wick2017variational}
T.~Wick.
\newblock Variational-monolithic ale fluid--structure interaction: Comparison
  of computational cost and mesh regularity using different mesh motion
  techniques.
\newblock In {\em Modeling, Simulation and Optimization of Complex Processes
  HPSC 2015}, pages 261--275. Springer, 2017.

\bibitem{richter2015monolithic}
T.~Richter.
\newblock A monolithic geometric multigrid solver for fluid--structure
  interactions in {ALE} formulation.
\newblock {\em International Journal for Numerical Methods in Engineering},
  2015.

\bibitem{burman2013explicit}
Erik Burman and Miguel{\hspace{0.167em}}A. Fern{\'{a}}ndez.
\newblock Explicit strategies for incompressible fluid-structure interaction
  problems: Nitsche type mortaring versus robin-robin coupling.
\newblock {\em International Journal for Numerical Methods in Engineering},
  97(10):739--758, December 2013.

\bibitem{causin2005added}
P.~Causin, J.-F. Gerbeau, and F.~Nobile.
\newblock {Added-mass effect in the design of partitioned algorithms for
  fluid--structure problems}.
\newblock {\em Computer methods in applied mechanics and engineering},
  194(42):4506--4527, 2005.

\bibitem{kuttler2008fixed}
U.~K{\"u}ttler and W.~Wall.
\newblock Fixed-point fluid--structure interaction solvers with dynamic
  relaxation.
\newblock {\em Computational mechanics}, 43(1):61--72, 2008.

\bibitem{landajuela2016}
M.~Landajuela, M.~Vidrascu, D.~Chapelle, and M.~A. Fern{\'{a}}ndez.
\newblock Coupling schemes for the {FSI} forward prediction challenge:
  Comparative study and validation.
\newblock {\em International Journal for Numerical Methods in Biomedical
  Engineering}, 33(4):e2813, July 2016.

\bibitem{gerbeau2005fluid}
J.-F. Gerbeau, M.~Vidrascu, and P.~Frey.
\newblock {Fluid--structure interaction in blood flows on geometries based on
  medical imaging}.
\newblock {\em Computers \& Structures}, 83(2):155--165, 2005.

\bibitem{wall2008fluid}
W.~A. Wall and T.~Rabczuk.
\newblock Fluid--structure interaction in lower airways of {CT}-based lung
  geometries.
\newblock {\em International Journal for Numerical Methods in Fluids},
  57(5):653--675, 2008.

\bibitem{jodlbauer2019parallel}
D~Jodlbauer, U~Langer, and T~Wick.
\newblock Parallel block-preconditioned monolithic solvers for fluid-structure
  interaction problems.
\newblock {\em International Journal for Numerical Methods in Engineering},
  117(6):623--643, 2019.

\bibitem{murea2017updated}
C.~Murea and S.~Sy.
\newblock Updated {L}agrangian/{A}rbitrary {L}agrangian--{E}ulerian framework
  for interaction between a compressible neo-{H}ookean structure and an
  incompressible fluid.
\newblock {\em International Journal for Numerical Methods in Engineering},
  109(8):1067--1084, 2017.

\bibitem{badia2008modular}
S.~Badia, A.~Quaini, and A.~Quarteroni.
\newblock Modular vs. non-modular preconditioners for fluid--structure systems
  with large added-mass effect.
\newblock {\em Computer Methods in Applied Mechanics and Engineering},
  197(49-50):4216--4232, 2008.

\bibitem{yang2015modeling}
K.~Yang, P.~Sun, L.~Wang, J.~Xu, and L.~Zhang.
\newblock Modeling and numerical studies for fluid--structure interaction
  involving an elastic rotor.
\newblock {\em Computer Methods in Applied Mechanics and Engineering}, 311:788
  -- 814, 2016.

\bibitem{xu2015well}
J.~Xu and K.~Yang.
\newblock Well-posedness and robust preconditioners for discretized
  fluid--structure interaction systems.
\newblock {\em Computer Methods in Applied Mechanics and Engineering},
  292:69--91, 2015.

\bibitem{holzapfel2001biomechanics}
G.~Holzapfel et~al.
\newblock Biomechanics of soft tissue.
\newblock {\em The handbook of materials behavior models}, 3:1049--1063, 2001.

\bibitem{OlshanskiiPetersReusken2006}
M.~A. Olshanskii, J.~Peters, and A.~Reusken.
\newblock Uniform preconditioners for a parameter dependent saddle point
  problem with application to generalized {S}tokes interface equations.
\newblock {\em Numerische Mathematik}, 105(1):159--191, 2006.

\bibitem{wichrowski2022}
M.~Wichrowski and P.~Krzyzanowski.
\newblock A matrix-free multilevel preconditioner for the generalized {S}tokes
  problem with discontinuous viscosity.
\newblock {\em Journal of Computational Science}, 63:101804, 2022.

\bibitem{jodlbauer2022matrixfree}
D.~Jodlbauer, U.~Langer, T.~Wick, and W.~Zulehner.
\newblock Matrix-free monolithic multigrid methods for {S}tokes and generalized
  {S}tokes problems, 2022.

\bibitem{Kronbichler2012}
M.~Kronbichler and K.~Kormann.
\newblock A generic interface for parallel cell-based finite element operator
  application.
\newblock {\em Computers \& Fluids}, 63:135--147, June 2012.

\bibitem{turek2006proposal}
S.~Turek and J.~Hron.
\newblock {\em Proposal for numerical benchmarking of fluid--structure
  interaction between an elastic object and laminar incompressible flow}.
\newblock Springer, 2006.

\bibitem{helenbrook2003mesh}
B.~Helenbrook.
\newblock Mesh deformation using the biharmonic operator.
\newblock {\em International journal for numerical methods in engineering},
  56(7):1007--1021, 2003.

\bibitem{ogden1972large}
R.~W. Ogden.
\newblock Large deformation isotropic elasticity -- on the correlation of
  theory and experiment for incompressible rubberlike solids.
\newblock {\em Proceedings of the Royal Society of London A: Mathematical,
  Physical and Engineering Sciences}, 326:565--584, 1972.

\bibitem{bigoni2012nonlinear}
D.~Bigoni.
\newblock {\em Nonlinear Solid Mechanics: Bifurcation Theory and Material
  Instability}.
\newblock Cambridge University Press, 2012.

\bibitem{cash1980integration}
J.R. Cash.
\newblock On the integration of stiff systems of odes using extended backward
  differentiation formulae.
\newblock {\em Numerische Mathematik}, 34(3):235--246, 1980.

\bibitem{heil2008solvers}
M.~Heil, A.~L. Hazel, and J.~Boyle.
\newblock {Solvers for large-displacement fluid--structure interaction
  problems: segregated versus monolithic approaches}.
\newblock {\em Computational Mechanics}, 43(1):91--101, 2008.

\bibitem{langer2018numerical}
U.~Langer and H.~Yang.
\newblock Numerical simulation of fluid--structure interaction problems with
  hyperelastic models: A monolithic approach.
\newblock {\em Mathematics and Computers in Simulation}, 145:186--208, 2018.

\bibitem{lozovskiy2015unconditionally}
A.~Lozovskiy, M.A. Olshanskii, V.~Salamatova, and Y.~V. Vassilevski.
\newblock {An unconditionally stable semi-implicit FSI finite element method}.
\newblock {\em Computer Methods in Applied Mechanics and Engineering},
  297:437--454, 2015.

\bibitem{van1986second}
J.J.I.M. van Kan.
\newblock {A second-order accurate pressure-correction scheme for viscous
  incompressible flow}.
\newblock {\em SIAM Journal on Scientific and Statistical Computing},
  7(3):870--891, 1986.

\bibitem{dong2010unconditionally}
S.~Dong and J.~Shen.
\newblock {An unconditionally stable rotational velocity-correction scheme for
  incompressible flows}.
\newblock {\em Journal of Computational Physics}, 229(19):7013--7029, 2010.

\bibitem{turek1996comparative}
S.~Turek.
\newblock A comparative study of time-stepping techniques for the
  incompressible {N}avier--{S}tokes equations: from fully implicit non-linear
  schemes to semi-implicit projection methods.
\newblock {\em International Journal for Numerical Methods in Fluids},
  22(10):987--1011, 1996.

\bibitem{boffi2013mixed}
D.~Boffi, F.~Brezzi, and M.~Fortin.
\newblock {\em Mixed finite element methods and applications}, volume~44.
\newblock Springer, 2013.

\bibitem{brooks1982streamline}
A.~N. Brooks and T.~Hughes.
\newblock Streamline upwind/{P}etrov--{G}alerkin formulations for convection
  dominated flows with particular emphasis on the incompressible
  {N}avier--{S}tokes equations.
\newblock {\em Computer methods in applied mechanics and engineering},
  32(1-3):199--259, 1982.

\bibitem{john2006discontinuity}
V.~John and P.~Knobloch.
\newblock On discontinuity-capturing methods for convection-diffusion
  equations.
\newblock In {\em Numerical mathematics and advanced applications}, pages
  336--344. Springer, 2006.

\bibitem{gee2011truly}
M.~W. Gee, U.~K{\"u}ttler, and W.~A. Wall.
\newblock {Truly monolithic algebraic multigrid for fluid--structure
  interaction}.
\newblock {\em International Journal for Numerical Methods in Engineering},
  85(8):987--1016, 2011.

\bibitem{wu2014fully}
Y.~Wu and X-Ch Cai.
\newblock A fully implicit domain decomposition based {ALE} framework for
  three-dimensional fluid--structure interaction with application in blood flow
  computation.
\newblock {\em Journal of Computational Physics}, 258:524--537, 2014.

\bibitem{langer2016recent}
U.~Langer and H.~Yang.
\newblock Recent development of robust monolithic fluid--structure interaction
  solvers.
\newblock {\em Fluid--Structure Interactions. Modeling, Adaptive Discretization
  and Solvers, Radon Series on Computational and Applied Mathematics}, 20,
  2016.

\bibitem{dealii2019design}
A.~Arndt, W.~Bangerth, D.~Davydov, T.~Heister, L.~Heltai, M.~Kronbichler,
  M.~Maier, Pelteret J.-P., Turcksin B., and Wells D.
\newblock The {deal.II} finite element library: Design, features, and insights.
\newblock {\em Computers \& Mathematics with Applications}, 81:407--422, 2021.

\bibitem{KronbichlerEtAl2022}
M.~Kronbichler, D.~Sashko, and P.~Munch.
\newblock Enhancing data locality of the conjugate gradient method for
  high-order matrix-free finite-element implementations.
\newblock {\em The International Journal of High Performance Computing
  Applications}, page 109434202211078, 2022.

\bibitem{jodlbauer2020fracture}
D.~Jodlbauer, U.~Langer, and T.~Wick.
\newblock Parallel matrix-free higher-order finite element solvers for
  phase-field fracture problems.
\newblock {\em Mathematical and Computational Applications}, 25(3):40, 2020.

\bibitem{kronbichler2012high}
M.~Kronbichler, T.~Heister, and W.~Bangerth.
\newblock {High accuracy mantle convection simulation through modern numerical
  methods}.
\newblock {\em Geophysical Journal International}, 191(1):12--29, 2012.

\bibitem{africa2023lifexcfd}
P.~Claudio Africa, I.~Fumagalli, M.~Bucelli, A.~Zingaro, L.~Dede', and
  A.~Quarteroni.
\newblock lifex-cfd: an open-source computational fluid dynamics solver for
  cardiovascular applications, 2023.

\bibitem{davydov2020matrix}
D.~Davydov, J.-P Pelteret, D.~Arndt, M.~Kronbichler, and P.~Steinmann.
\newblock A matrix-free approach for finite-strain hyperelastic problems using
  geometric multigrid.
\newblock {\em International Journal for Numerical Methods in Engineering},
  121(13):2874--2895, 2020.

\bibitem{zulehner2000class}
W.~Zulehner.
\newblock A class of smoothers for saddle point problems.
\newblock {\em Computing}, 65(3):227--246, 2000.

\bibitem{ADAMS2003593}
M.~Adams, M.~Brezina, J.~Hu, and Tuminaro R.
\newblock Parallel multigrid smoothing: polynomial versus {G}auss--{S}eidel.
\newblock {\em Journal of Computational Physics}, 188(2):593 -- 610, 2003.

\end{thebibliography}

\end{document}